\documentclass[letterpaper,12pt]{amsart}
\usepackage{pict2e,amsfonts,color,latexsym,mathrsfs,amssymb,amsmath,epsfig,rotating}
\newtheorem{thm}{Theorem}[section]

\newtheorem{pellet}{Pellet's Theorem} 
 
\newtheorem{kapra}{\scalebox{.9}[1]{Kapranov's Non-Archimedean Amoeba Theorem}} 
 
\newtheorem{dfn}[thm]{Definition}

\newtheorem{cor}[thm]{Corollary} 
\newtheorem{rem}[thm]{Remark}
\newtheorem{prop}[thm]{Proposition}
\newtheorem{lemma}[thm]{Lemma}
\newtheorem{ex}[thm]{Example}
\newlength{\smale}
\newlength{\jmr}
\newlength{\khov}
\newlength{\bernd}
\makeatletter
\renewcommand*\env@matrix[1][c]{\hskip -\arraycolsep
  \let\@ifnextchar\new@ifnextchar
  \array{*\c@MaxMatrixCols #1}}
\makeatother

\definecolor{red}{rgb}{.5,0,0} 
\definecolor{green}{rgb}{0,.4,0} 
\definecolor{blue}{rgb}{0,0,.5} 
\newcommand{\vm}{v_{\mathrm{min}}}
\newcommand{\vp}{v_{\mathrm{max}}}
\newcommand{\floor}[1]{\left\lfloor#1\right\rfloor} 
\newcommand{\ceil}[1]{\left\lceil#1\right\rceil} 
\newcommand{\thth}{^{\text{\underline{th}}}}

\newcommand{\np}{{\mathbf{NP}}}

\newcommand{\pp}{\mathbf{P}}

\newcommand{\pspa}{\mathbf{PSPACE}}
\newcommand{\expt}{\mathbf{EXPTIME}}
\newcommand{\eps}{\varepsilon}
\newcommand{\Pro}{{\mathbb{P}}}
\newcommand{\Log}{\mathrm{Log}}

\newcommand{\Q}{\mathbb{Q}}
\newcommand{\R}{\mathbb{R}}

\newcommand{\C}{\mathbb{C}}
\newcommand{\N}{\mathbb{N}}
\newcommand{\Z}{\mathbb{Z}}
\newcommand{\bO}{\mathbf{O}}

\newcommand{\Zn}{\Z^n}

\newcommand{\Qn}{\Q^n}
\newcommand{\Rn}{\R^n}

\newcommand{\Cn}{\C^n}

\newcommand{\Cs}{\C^*}

\newcommand{\Csn}{{(\C^*)}^n}

\renewcommand{\qed}{$\blacksquare$}
\newcommand{\dia}{$\diamond$}
\newcommand{\amoeba}{\mathrm{Amoeba}}
\newcommand{\newt}{\mathrm{Newt}}
\newcommand{\anewt}{\mathrm{ArchNewt}}
\newcommand{\trop}{\mathrm{Trop}}
\newcommand{\atrop}{\mathrm{ArchTrop}}
\newcommand{\supp}{\mathrm{Supp}}
\newcommand{\conv}{\mathrm{Conv}}
\newcommand{\size}{\mathrm{size}}


\author{Mart\'in Avenda\~no}
\address{Department of Mathematics, University of Zaragoza, Zaragoza, 
Spain.}  
\email{mavendar@yahoo.com.ar}
\author{Roman Kogan}
\address{Department of Mathematics, Texas A\&M University TAMU 3368,
College Station, Texas \ 77843-3368, USA.} 
\email{romwell@math.tamu.edu}
\author{Mounir Nisse}
\address{Department of Mathematics,
Texas A\&M University TAMU 3368,
College Station, Texas \ 77843-3368, USA.}  
\email{nisse@math.tamu.edu}
\author{J.\ Maurice Rojas}
\address{Department of Mathematics,
Texas A\&M University TAMU 3368,
College Station, Texas \ 77843-3368, USA.}  
\email{rojas@math.tamu.edu}
\thanks{Partially supported by NSF grants DMS-0915245, 
CCF-1409020, and DMS-1460766. 
J.M.R.\ was also partially supported by 
DOE ASCR grant DE-SC0002505 and Sandia National Laboratories. 
Part of this work was presented as a talk at MEGA 2013 
(June 4--7, Goethe University, Frankfurt, Germany). }

\title[Archimedean Amoebae and Tropical Hypersurfaces]{\mbox{}\\
\vspace{-1.25in}
Metric Estimates and Membership Complexity for Archimedean Amoebae and 
Tropical Hypersurfaces}  

\date{\today} 

\begin{document}

\begin{abstract}  
Given any complex Laurent polynomial $f$, 
$\amoeba(f)$ is the image of its\linebreak complex zero set under the
coordinate-wise log absolute value map. We give an efficiently constructible
polyhedral approximation, $\atrop(f)$, of $\amoeba(f)$, and derive 
explicit upper and lower bounds, solely as a function 
of the number of monomial terms of $f$, for the Hausdorff distance between 
these two sets. We also show that deciding whether a given point lies in 
$\atrop(f)$ is doable in polynomial-time, for any fixed dimension, unlike the 
corresponding problem for $\amoeba(f)$, which is $\np$-hard already in one 
variable. $\atrop(f)$ can thus serve as a canonical low order 
approximation to start any higher order iterative polynomial system 
solving algorithm, such as homotopy continuation. $\atrop(f)$ also provides an 
Archimedean analogue of Kapranov's Non-Archimedean Amoeba Theorem  
and a higher-dimensional extension of earlier estimates of Mikhalkin and 
Ostrowski.  
\end{abstract} 

\maketitle

\vspace{-.4cm} 
\mbox{}\hfill {\bf In memory of Mikael Passare.} \hfill \mbox{}\\

\vspace{-.6cm} 
\section{Introduction}  

\noindent 
\scalebox{.944}[1]{One of the happiest coincidences in algebraic geometry is 
that the norms of roots of \mbox{polynomials}}\linebreak can be estimated 
through polyhedral geometry. Perhaps the earliest incarnation of this fact 
was Isaac Newton's use of a polygon to determine initial exponents of series 
expansions for\linebreak algebraic functions in one variable. This was 
detailed in a letter, dated October 24, 1676 \cite{newton}, that Newton wrote 
to Henry Oldenburg.  In modern terminology, Newton counted,\linebreak 
\scalebox{.94}[1]{with multiplicity, the $s$-adic valuations of roots of 
univariate polynomials over the {\em Puiseux series}}\linebreak 
{\em field} $\C\langle\langle s\rangle \rangle$ (see, e.g., Theorem 
\ref{thm:newton} from Section \ref{sub:arch} below). 
Newton's result has since been extended to arbitrary non-Archimedean 
fields (see, e.g., \cite{dumas,weiss}). Tropical\linebreak   
geometry (see, e.g., \cite{bergman,ekl,litvinov,tropical1,bakerrumely,barba,
macsturmf}) continues to deepen the links between algebraic, arithmetic, and 
polyhedral geometry, but has so far concentrated mainly on algebraic sets 
over fields other than $\C$. 

We will use tropical methods to efficiently approximate 
complex algebraic hypersurfaces in arbitrary dimension 
(see Theorem \ref{thm:finally} and Corollary \ref{cor:mixed} 
in Sections \ref{sub:n}--\ref{sub:apps} below), and derive an 
Archimedean analogue of Newton's result along the way (Theorem 
\ref{thm:harder} in Section \ref{sub:1} below).  
While our approximations can be coarse, their computational cost is 
quite low (see Theorem \ref{thm:cxity2} in Section \ref{sub:cxity} below), 
and initial experiments indicate that they are often good enough to yield 
high-quality start points for homotopy algorithms applied to sparse polynomial 
systems (see, e.g., \cite{aggr}).   
\begin{dfn}
\label{dfn:archtrop}
Let $\Cs\!:=\!\C\setminus\{0\}$, let $c_1,\ldots,c_t\!\in\Cs$, and call any
$f\!\in\!\C\!\left[x^{\pm 1}_1,\ldots,x^{\pm 1}_n\right]$ of the\linebreak 
form $f(x)\!=\!\sum^t_{i=1} c_i x^{a_i}$, 
with $\{a_1,\ldots,a_t\}$ of cardinality 
$t\!\geq\!1$, an {\em $n$-variate $t$-nomial}. (The\linebreak notation 
$x\!=\!(x_1,\ldots,x_n)$
and $x^{a_i}\!=\!x^{a_{1,i}}_1\cdots x^{a_{n,i}}_n$ is understood.)
We then define the\linebreak 
{\em (ordinary) Newton polytope} of $f$ to be
$\newt(f)\!:=\!\conv\!\left(\{a_i\}_{i\in [t]}\right)$, and the {\em
Archimedean} Newton polytope of $f$ to be
$\anewt(f)\!:=\!\conv\!\left(\{(a_i,-\log|c_i|)\}_{i\in [t]}\right)$. \dia 
\end{dfn}

\noindent 
Hadamard defined the $n\!=\!1$ case of $\anewt(f)$ around 1893, and observed 
a\linebreak  
relationship between the absolute values of the complex roots of $f$ 
and the slopes of certain edges of $\anewt(f)$ 
\cite[pp.\ 174--175 \& 201]{hadamard} (see also \cite[pp.\ 120--121]{ostrowski} 
and \cite[Ch.\ IX, pp.\ 193--202]{valiron}). We'll see below that, for 
$n\!\geq\!2$, approximating absolute\linebreak  
values of complex roots can be reduced to maximizing certain linear 
forms over $\anewt(f)$, and this ultimately leads us to a particular class of 
tropical varieties. In what follows, for any 
$\zeta\!=\!(\zeta_1,\ldots,\zeta_n)\!\in\!\Csn$, we set 
$\Log|\zeta|\!:=\!(\log|\zeta_1|,\ldots,\log|\zeta_n|)$ and define 
$\amoeba(f):=\left\{\Log|\zeta|\; | \; f(\zeta)\!=\!0 \text{ and } 
\zeta\!\in\!\Csn\right\}$. Clearly, $\amoeba(f)$ is empty 
when $f$ is a monomial, so we assume $t\!\geq\!2$ henceforth. 
\begin{ex}
\label{ex:binom}
For the trivariate binomial $f(x)\!=\!7-x^2_1x^3_2x^5_3$, it 
is easily checked that $\anewt(f)$ is the line segment in $\R^4$ 
connecting $(0,0,0,-\log 7)$ and $(2,3,5,0)$. In\linebreak 
particular, $\Log|\zeta|\!\in\!\amoeba(f) \Longleftrightarrow 
7\!=\!\left|\zeta^2_1 \zeta^3_2 \zeta^5_3\right|$, and 
thus it is clear that $\amoeba(f)$ is exactly the affine hyperplane in 
$\R^3$ defined by $2v_1+3v_2+5v_3=\log 7$. Note also that any 
$(v_1,v_2,v_3)\!\in\!\amoeba(f)$ makes the vector $(v_1,v_2,v_3,-1)$ 
perpendicular to $\anewt(f)$. \dia  
\end{ex}

\begin{ex}
\label{ex:tri} 
When $f(x_1)\!:=\!\frac{1}{89}-x^{16}_1+x^{49}_1$ it turns out that 
$\amoeba(f)$ consists of exactly $26$\linebreak 
\scalebox{.945}[1]{points.  
However, the points of $\amoeba(f)$ cluster tightly 
about just $\pmb{2}$ values: Exactly 
$16$ complex}\linebreak 
\scalebox{.95}[1]{roots of $f$ have norm near 
$\sqrt[16]{\frac{1}{89}}\!\approx\!0.7553...$ (to at 
least $4$ decimal places) and exactly $33$ complex}\linebreak 
roots of $f$ have norm near 
$1$ (to $3$ decimal places). Here, $\anewt(f)$ is the convex hull of 
$\left\{\left(0,-\log\frac{1}{89}\right),(16,0),(49,0)
\right\}$, which is the triangle drawn below. Note also that the 
linear\linebreak  

\vspace{-1.65cm}
\noindent
\begin{minipage}[b]{0.5\linewidth}
\vspace{0pt}
\scalebox{.52}[.52]{\begin{picture}(500,145)(0,-20)
\put(0,44.886363){{\color{black}\circle*{3}}}
\put(0,44.886363){\line(16,-4.4886){160}}
\put(160,0){{\color{black}\circle*{3}}}
\put(160,0){\line(1,0){330}}
\put(490,0){{\color{black}\circle*{3}}}
\put(490,0){\line(-49,4.486){490}}
\end{picture}} 
\end{minipage}\hspace{.7cm}
\begin{minipage}[b]{0.45\linewidth}
\vspace{0pt} 
form $\textcolor{blue}{\frac{\log\frac{1}{89}}{16}}v_1-v_2$ is 
maximized on the lower left edge of $\anewt(f)$, while the linear 
form $\textcolor{blue}{0}v_1-v_2$ is maximized on the 
\end{minipage}

\noindent
lower right edge of $\anewt(f)$ (if we restrict both linear forms to 
$\anewt(f)$). More than coincidentally, every point of $\amoeba(f)$ 
is within $0.00034$ of some point of $\left\{\textcolor{blue}{\frac{1}{16}
\log\frac{1}{89}},\textcolor{blue}{0}\right\}$, and the horizontal lengths 
($16$ and $33$) of the two lower edges 
count the number of roots with norm in the corresponding cluster.
Note also that when $\log|\zeta|\!=\!\frac{1}{16}\log\frac{1}{89}$ we have 
$\frac{1}{89}\!=\!\left|-\zeta^{16}\right|\!>\!\left|\zeta^{49}
\right|$, and when $\log|\zeta|\!=\!0$ have $\frac{1}{89}\!<\!
\left|-\zeta^{16}\right|\!=\!\left|\zeta^{49}\right|$. 
Furthermore, when $\log|\zeta|\!\not\in\!\left\{\frac{1}{16}\log\frac{1}{89},
0\right\}$, the set $\left\{\left|\frac{1}{89}\right|,\left|-\zeta^{16}
\right|,\left|\zeta^{49}\right|\right\}$ has cardinality $3$. \dia
\end{ex}

We refer the reader to the outstanding texts \cite{ziegler,triang} for 
further background on polytopes, faces, and normal fans. 
\begin{dfn} 
We define the {\em Archimedean} tropical variety of $f$,
$\atrop(f)$, to be\\  
\mbox{}\hfill $\left\{v\!\in\!\Rn \; | \; (v,-1) \text{ {\em is an 
outer normal of a positive-dimensional face of} } \anewt(f) 
\right\}$\hfill\mbox{}\\ 
when $t\!\geq\!2$ and, when $t\!=\!1$, we set $\atrop(f)\!=\!\emptyset$. 
We also call a face of $\anewt(f)$\\ 
a {\em lower} face iff it has an outer normal of the form $(v,-1)$ for some 
$v\!\in\!\Rn$. \dia  
\end{dfn}

\noindent 
For instance, in Examples \ref{ex:binom} and \ref{ex:tri}, $\atrop(f)$ 
was, respectively, a plane in $\R^3$ equal to $\amoeba(f)$, and then 
a pair of points around which $\amoeba(f)$ clustered. 
While $\atrop(f)$ has appeared under different guises in earlier work 
(see, e.g., \cite{ostrowski,mik3,passarerullgard,prs,theobalddewolff}), 
explicit metric estimates for how well $\atrop(f)$ approximates $\amoeba(f)$ in 
{\em arbitrary} dimension have not yet appeared in the literature.\footnote{The 
only works in this direction that we are aware of are \cite[Rem.\ 4.2, 
Inequality 25]{ags} (for $n\!=\!1$, using new bounds for matrix polynomials)  
and \cite[Lemma 8.5, Pg.\ 360]{mikhalkin} (for the case $n=2$).} 

\subsection{Metric Estimates, With Multiplicity, in One Variable} 
\label{sub:1} We now give explicit bounds on how $\amoeba(f)$ 
clusters about the points of $\atrop(f)$. 
In what follows, for any line segment $L\!\subset\!\R^2$ with vertices $(a,b)$ 
and $(c,d)$, its {\em horizontal length} is $\lambda(L)\!:=\!|c-a|$.  
\begin{thm} 
\label{thm:harder} 
Given any univariate $t$-nomial $f$ with $t\!\geq\!2$, let 
$\alpha_f\!:=\!\min\atrop(f)-\log 2$, 
$\beta_f\!:=\!\max\atrop(f)+\log 2$, let  
$\Gamma$ be any connected component of the union of open intervals 
$U_f\!:=\!(\alpha_f,\beta_f) \cap \bigcup\limits_{v\in \atrop(f)} 
\!\!\!\!\!\!\! 
(v-\log 3,v+\log 3)$, and let $\Lambda_\Gamma$ be the sum of $\lambda(L)$ over 
all edges $L$ of $\anewt(f)$ with outer normal $(v,-1)$ satisfying 
$v\!\in\!\Gamma$. Then the number of roots $\zeta\!\in\!\C$ of $f$ with 
$\log|\zeta|\!\in\!\Gamma$, counting multiplicity, is exactly 
$\Lambda_\Gamma$. In particular, $\amoeba(f)\!\subset\!U_f$ and 
$\Lambda_\Gamma\!\geq\!1$.  
\end{thm} 

\noindent 
Theorem \ref{thm:harder} is proved in Section \ref{sub:harder}, 
where a slight sharpening is also provided for $t\!=\!3$. 
\begin{ex}
If $f(x_1)\!:=\!1+19162399831x^{16}_1+x^{49}_1$ then $U_f$ is a disjoint 
union of two intervals, and Theorem \ref{thm:harder} 
tells us that $f$ has exactly $16$ (resp.\ $33$) roots with log-norm in the 
open interval $(-2.17292,-0.381151202192)$ (resp.\  
$(-0.381151202190,1.41061)$). In fact, for this example, the much 
smaller sub-intervals\\ 
\mbox{}\hfill $\frac{-\log 19162399831}{16}+10^{-32}(-1,1)$ \ \ \ and \ \ \ 
$\frac{\log(19162399831)}{33}+10^{-16}(-1,1)$ \hfill\mbox{}\\ 
(respectively centered at 
the $2$ points of $\atrop(f)$) still respectively contain the 
same number of log-norms. \dia 
\end{ex} 
\begin{ex} 
The constants in the definition of $U_f$ from Theorem 
\ref{thm:harder} are optimal: Assertion (c) of Corollary 
\ref{cor:extremecauchy} in Section \ref{sec:back} below reveals 
that the $\log 2$ in the definition of $\alpha_f$ and $\beta_f$ can 
not be replaced by any smaller constant, and Lemma 
\ref{lemma:opt} 
from Section \ref{sec:back} shows that the $\log 3$ can not be replaced by 
any smaller constant. \dia 
\end{ex} 

In Section \ref{sec:back} we 
also discuss how the neighborhood $U_f$ improves (or complements) 
earlier root norm estimates in \cite{hadamard,ostrowski,ags}, and how the 
$\Lambda_\Gamma$ provide an {\em Archimedean} version of tropical intersection 
multiplicity. 

\subsection{Explicit Piecewise Linear Approximations of Hypersurfaces}  
\label{sub:n} 
Moving on to the multivariate case, let us first make some 
basic observations on the structure of $\atrop(f)$. 
\begin{prop} 
\label{prop:dumb} 
If $f$ is an $n$-variate binomial then $\amoeba(f)$ and $\atrop(f)$ 
are\linebreak 
identical affine hyperplanes in $\Rn$. \qed 
\end{prop} 
\begin{lemma} 
\label{lemma:basic} 
Suppose $f$ is an $n$-variate $t$-nomial with $\newt(f)$ of dimension $k$. 
Then:\\
\mbox{}\hspace{.5cm} (0) $k\!\leq\!\min\{n,t-1\}$.\\ 
\mbox{}\hspace{.5cm} (1) $k\!=\!1 \Longrightarrow \atrop(f)$ is a non-empty 
disjoint union of at most $t-1$ parallel affine\\ 
\mbox{}\hspace{1.1cm}hyperplanes in $\Rn$.\\
\mbox{}\hspace{.5cm} (2) $k\!\geq\!2 \Longrightarrow \atrop(f)$ is a 
path-connected polyhedral complex, of pure dimension\\ 
\mbox{}\hspace{1.3cm}$n-1$, with at most $t(t-1)/2$ faces of dimension 
$n-1$.\\  
\mbox{}\hspace{.5cm} (3) $t\!=\!k+1 \Longrightarrow 
\atrop(f)\!\subseteq\!\amoeba(f)$ 
and both $\amoeba(f)$ and $\atrop(f)$ are\\ 
\mbox{}\hspace{1.3cm}contractible. \qed 
\end{lemma} 

\noindent 
Proposition \ref{prop:dumb} is elementary. Assertions (0)--(2) 
of Lemma \ref{lemma:basic} follows easily from the definition 
of $\atrop(f)$, thanks to polyhedral duality \cite{ziegler}. 
(See also \cite[Ch.\ 3, Sec.\ 3]{macsturmf} for a much more detailed 
discussion in the non-Archimedean setting.) Assertion (3) of Lemma 
\ref{lemma:basic} was one of the first basic 
topological results on amoebae and can be found, for instance, in 
\cite[Prop.\ 3.1.8]{forsbergthesis}, \cite[Thms.\ 8 \& 12]{rullgardthesis}, 
and \cite[Lemma 3.4 (a)]{theobalddewolff}.  

Our next main result is that every point of $\amoeba(f)$ is within 
an explicit distance of some point of $\atrop(f)$, {\em and vice-versa}, 
independent of the degree or number of variables of $f$. We use $|\cdot|$ for 
the standard $\ell_2$-norm on $\Cn$. 
\begin{dfn} 
For any $\eps\!>\!0$ and $X\!\subseteq\!\Rn$ we define the {\em open 
$\eps$-neighborhood of $X$} to be $X_\eps\!:=\!\{x\!\in\!\Rn\; | \; 
|x-x'|\!<\!\eps \text{ for some } x'\!\in\!X\}$, and let 
$\overline{X}_\eps$ denote its Euclidean closure. \dia 
\end{dfn} 
\begin{thm}
\label{thm:finally} 
For any $f\!\in\!\C\!\left[x^{\pm 1}_1,\ldots,x^{\pm 1}_n\right]$ with exactly 
$t\!\geq\!2$ monomial terms and $\newt(f)$ of dimension $k$ we have:  
\begin{enumerate} 
\item{(a) $\amoeba(f)\subseteq\overline{\atrop(f)}_{\log(t-1)}$ and, 
for $k\!=\!1$, $\amoeba(f)\subsetneqq\atrop(f)_{\log 3}$.\\  
(b) $\atrop(f)\subseteq\amoeba(f)_{\eps_{k,t}}$ where 
$\eps_{1,t}\!:=\!(\log 9)t -\log\frac{81}{2}\!<\!2.2t-3.7$,\\  
$\eps_{2,t}\!:=\!\sqrt{2}(t-2)\left((\log 9)t-\log\frac{81}{2}\right)\! 
<\!(t-2)(3.11t-5.23)$, and\\ 
$\eps_{k,t}\!:=\!\sqrt{k}\ceil{\frac{1}{4}
t(t-1)}\left((\log 9)t-\log\frac{81}{2}\right)$ for $k\!\geq\!3$.}  
\item{Let $\varphi(x)\!:=\!1+x_1+\cdots+x_{t-1}$ and 
$\psi(x)\!:=\!(x_1+1)^{t-k} +x_2+\cdots+x_k$. Then\\ 
(a) $\amoeba(\varphi)$ contains a point at 
distance $\log(t-1)$ from $\atrop(\varphi)$ and\\  
(b) $\atrop(\psi)$ contains points approaching distance 
$\log(t-k)$ from $\amoeba(\psi)$.} 
\end{enumerate} 
\end{thm}
\noindent 
We prove Theorem \ref{thm:finally} in Section \ref{sec:finally}. 
For multivariate polynomials, our bounds appear to be the first allowing 
dependence on just the number of terms $t$. 
In particular, Assertion (1a) sharpens, and extends to arbitrary 
dimension, an earlier bound of Mikhalkin for the case $n\!=\!2$: 
Letting $L$ denote the number of lattice points in the Newton polygon of $f$, 
\cite[Lemma 8.5, Pg.\ 360]{mikhalkin} asserts that $\amoeba(f)$ is contained 
in the possibly larger neighborhood $\overline{\atrop(f)}_{\log(L-1)}$.  
Assertion (2a) of Theorem \ref{thm:finally} shows that  
the size of the neighborhood from Assertion (1a) is in fact optimal 
for the infinite family of cases $t\!=\!k+1\!\geq\!3$.  

Finding the tightest neighborhood of $\amoeba(f)$ containing 
$\atrop(f)$ appears to be an open problem: We are unaware of any 
earlier multivariate version of Assertion (1b). The only other earlier 
distance bound between an amoeba (of positive dimension) and a polyhedral 
approximation we know 
of is a result of Viro \cite[Sec.\ 1.5]{virologpaper} on the distance between 
the graph of a univariate polynomial (drawn on log paper) and a piecewise 
linear curve that is ultimately a piece of the $n\!=\!2$ case of $\atrop(f)$ 
here.

\vspace{-.4cm} 
\noindent 
\begin{picture}(200,200)(5,-105)
\put(385,-105){
\begin{minipage}[b]{0.3\linewidth}
\vspace{0pt}
\epsfig{file=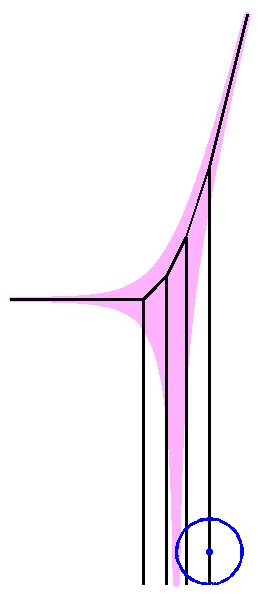,height=2.7in,clip=}
\end{minipage}} 
\put(0,-100){
\begin{minipage}[b]{.9\linewidth}
\vspace{0pt}
\begin{ex}
\label{ex:archoff}  
\scalebox{.97}[1]{Setting $\psi(x)\!=\!(x_1+1)^4+x_2$ we see $\amoeba(\psi)
\cap([-7,7]\times [-12,12])$}\linebreak 
\scalebox{.98}[1]{and \ $\atrop(\psi)\cap([-7,7]\times [-12,12])$ \ on the 
right.  \ $\atrop(\psi)$ contains the ray}\linebreak 
\scalebox{.98}[1]{$(\log 4,4\log 4)+\R_+(0,-1)$ and this rightmost 
downward-pointing ray contains points}\linebreak 
\scalebox{.96}[1]{with distance from $\amoeba(\psi)$ 
approaching $\log 4$. We also observe that Viro's earlier}\linebreak  
\scalebox{.98}[1]{polygonal approximation of graphs of 
univariate polynomials on log paper, applied}\linebreak 
\scalebox{.98}[1]{here, would result in the polygonal curve that is the   
subcomplex of $\atrop(\psi)$}\linebreak
obtained by deleting all $4$ downward-pointing rays. \dia  
\end{ex} 

\vspace{-.1cm} 
\mbox{}\hspace{.5cm}It is worth comparing Theorem \ref{thm:finally} to two 
other methods for approximating complex amoebae: Purbhoo, in \cite{purbhoo}, 
describes a uniformly convergent sequence\\ 
\scalebox{1}[1]{of outer polyhedral approximations to any
amoeba, using cyclic resultants.$^2$ While}\linebreak 
$\atrop(f)$ lacks this refinability, the computation of $\atrop(f)$ is 
considerably simpler: See Section \ref{sub:cxity} below and \cite{aggr}. 
$\atrop(f)$ is actually closer in spirit to the {\em spine} of 
$\amoeba(f)$. The latter construction, based on a multivariate 
\end{minipage}}  
\end{picture} 
\addtocounter{footnote}{1} 
\footnotetext{See also \cite{timo} for recent computational 
improvements to Purbhoo's outer approximation.} 

\vspace{-.15cm} 
\noindent 
version of Jensen's Formula from complex analysis, is due to Passare and 
Rullg\aa{}rd \cite[Sec.\ 3]{passarerullgard} and results in a 
polyhedral complex that is always contained in, and is homotopy equivalent to,
$\amoeba(f)$. 
Unfortunately, the computational complexity of the spine is not as 
straightforward as that of $\atrop(f)$. 
Further background on the computational complexity of amoebae can
be found in \cite{theobald,boundewolf,wolfsos}. 

Our final main results concern the complexity of deciding whether 
a given point lies in a given amoeba or Archimedean tropical 
variety. However, let first us observe a consequence of our metric estimates 
for systems of polynomials. 
\subsection{Coarse, but Fast, Isolation of Roots of Polynomial Systems} 
\label{sub:apps} 
An immediate consequence of Assertion (1a) of Theorem \ref{thm:finally} is an  
estimate for isolating the possible norm vectors of complex roots of arbitrary 
systems of multivariate polynomial equations. 
\begin{cor} 
\label{cor:mixed} 
Suppose $f_1,\ldots,f_m\!\in\!\C\!\left[x^{\pm 1}_1,\ldots,x^{\pm 1}_n\right]$ 
where $f_i$ has exactly $t_i$ monomial terms for all $i$. Then any root 
$\zeta\!\in\!\Csn$ of $F\!=\!(f_1,\ldots,f_m)$ satisfies\\ 
\mbox{}\hfill $\Log|\zeta|\!\in\!\overline{\atrop(f_1)}_{\eps_1}  
\cap \cdots \cap \overline{\atrop(f_m)}_{\eps_m}$, \hfill\mbox{}\\ 
where $\eps_i\!:=\!\log(t_i-1)$ for all $i$. \qed 
\end{cor} 

\begin{ex} 
We can isolate the log-norm vectors 
of the complex roots of the $3\times 3$ system\\ 
\mbox{}\hfill  
$F\!:=\!(\textcolor{red}{f_1},\textcolor{green}{f_2},\textcolor{blue}{f_3}
)\!:=\!(\textcolor{red}{x_1x_2-x^2_1-1/16^6},\textcolor{green}{x_2x_3-1
-x^2_1/16^6},\textcolor{blue}{x_3-1-x^2_1/16^{18}})$ \hfill \mbox{}\\  
\scalebox{.92}[1]{via Corollary \ref{cor:mixed} as follows: Find the 
points of $X\!:=\!\textcolor{red}{\atrop(f_1)}\cap\textcolor{green}{
\atrop(f_2)}\cap\textcolor{blue}{\atrop(f_3)}$}\linebreak 
by searching through 
suitable triplets of edges of the $\anewt(f_i)$, and then create 
isolating parallelepipeds about the points of $X$. More precisely, 
observe that\\  
\mbox{}\hfill 
$\textcolor{red}{\conv(\{(1,1,0,0),(2,0,0,0)\})}$, 
$\textcolor{green}{\conv(\{(0,1,1,0),(0,0,0,0)\})}$, 
$\textcolor{blue}{\conv(\{(0,0,1,0),(0,0,0,0)\})}$\hfill\mbox{}\\ 
are respective edges of $\textcolor{red}{\anewt(f_1)}$, 
$\textcolor{green}{\anewt(f_2)}$, and $\textcolor{blue}{\anewt(f_3)}$, 
and the vector $(0,0,0,-1)$ is an outer normal to each of these edges. 
So $(0,0,0)$ is a point of $X$. Running through the remaining triplets 
we then obtain that $X$ in fact consists of exactly $4$ points:\\ 
\mbox{}\hfill $\Log\left|\left(\frac{1}{16^6},1,1\right)
\right|$ \ , \  $\Log|(1,1,1)|$ \ , \ $\Log\left|\left(16^6,16^6,1
\right)\right|$ \ , \ and \  
$\Log\left|\left(16^{12},16^{12},16^6\right)\right|$.\hfill \mbox{}\\  
So Corollary \ref{cor:mixed} tells us that the points of $Y\!:=\!
\textcolor{red}{\amoeba(f_1)} 
\cap\textcolor{green}{\amoeba(f_2)}\cap\textcolor{blue}{\amoeba(f_3)}$ lie in 
the union of the $4$ parallelepipeds drawn below to the right: Truncations of 
$\textcolor{red}
{\atrop(f_1)}$, $\textcolor{green}{\atrop(f_2)}$, and $\textcolor{blue}
{\atrop(f_3)}$ are drawn below on the left, and the middle illustration 
uses transparency to further detail the intersection. \\  
\mbox{}\hfill 
\epsfig{file=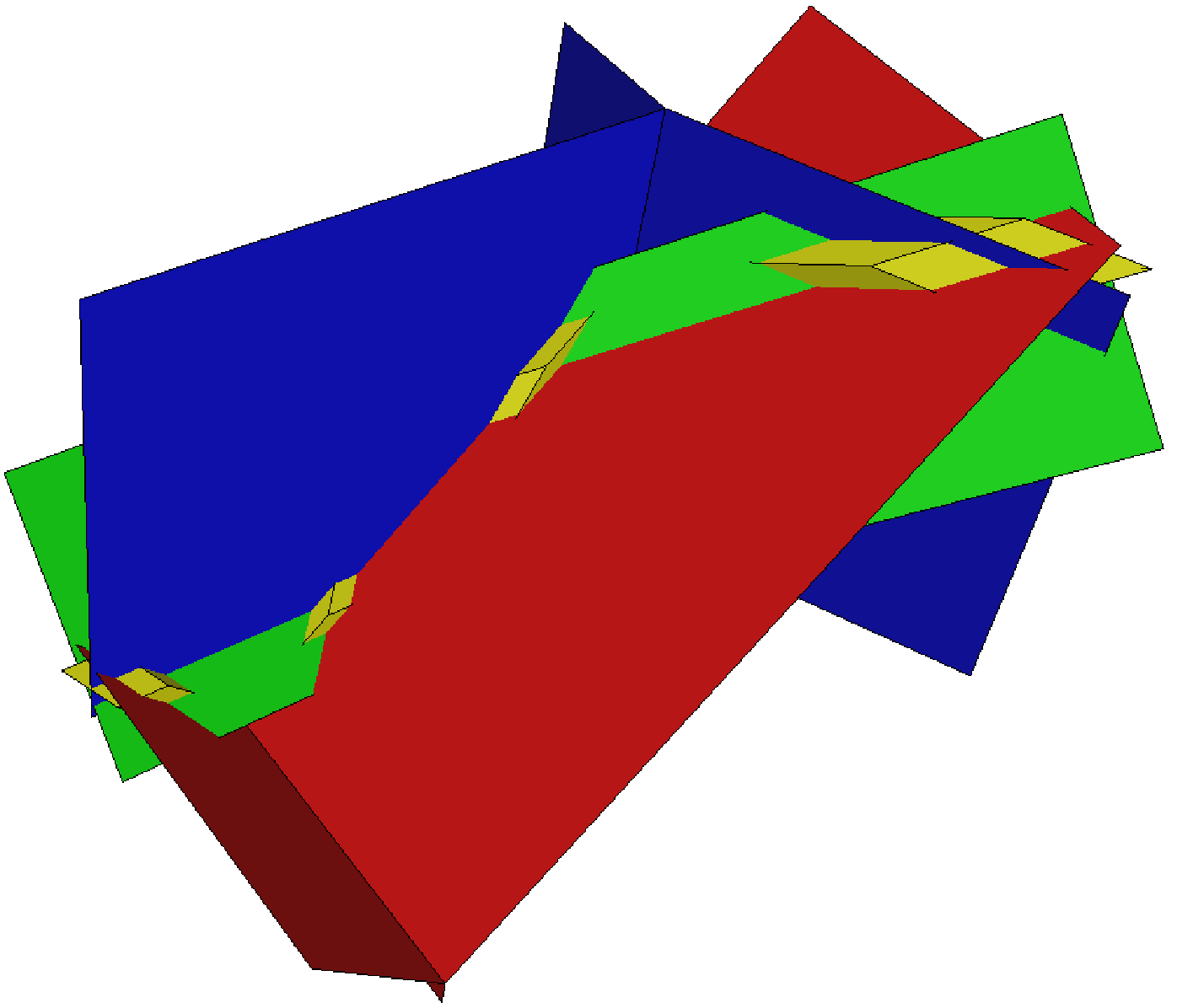,height=1.5in,clip=} 
\hfill 
\epsfig{file=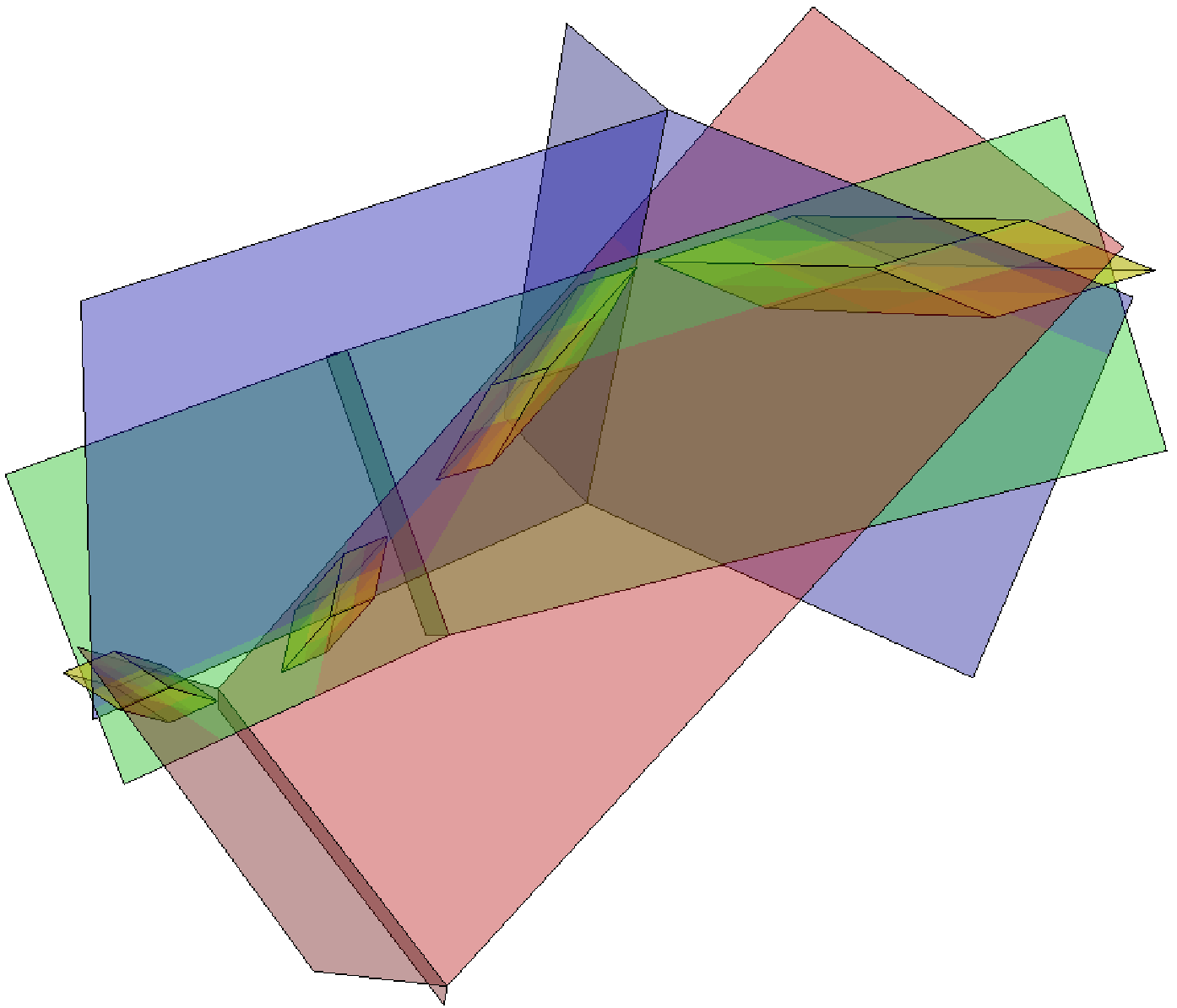,height=1.5in,clip=}
\hfill 
\raisebox{1cm}{\epsfig{file=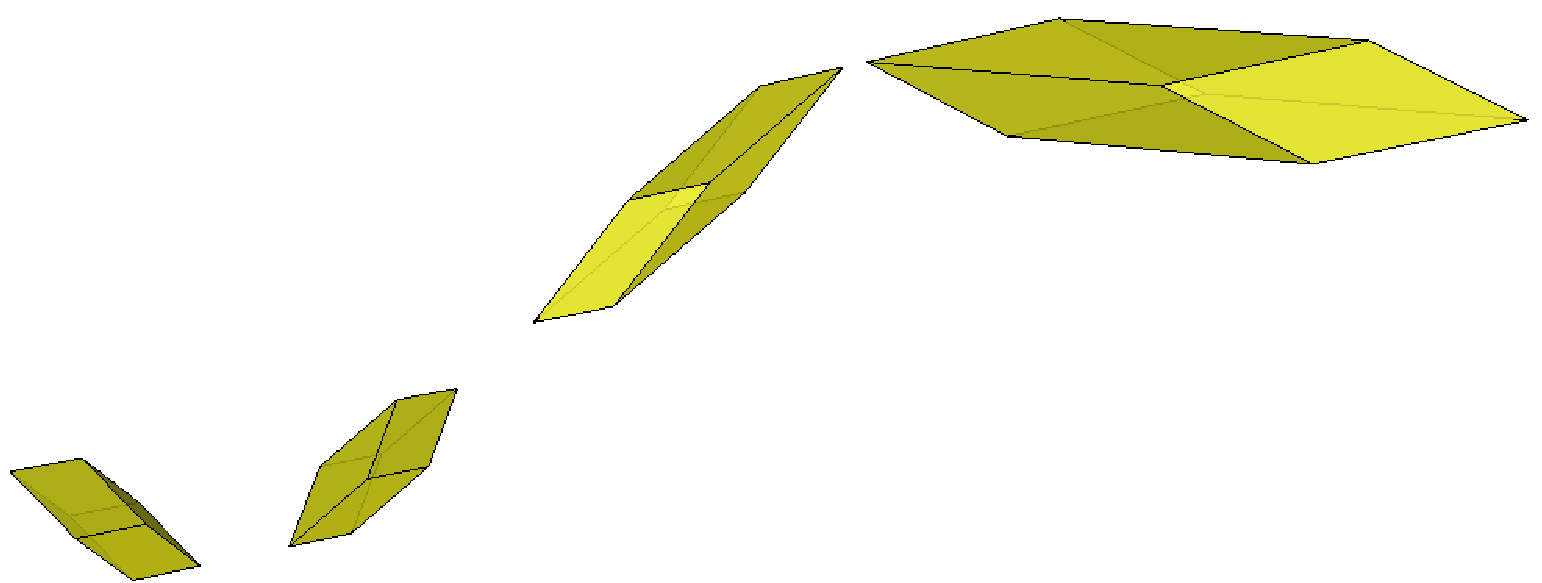,height=.8in,clip=}}
\hfill\mbox{}\\
Suitably ordered, each point of $X$ is actually within distance 
$0.11\times 10^{-6}$ ($<\!0.693...\!=\!\log 2$) of some point of 
$Y$ (and vice-versa), well in accordance with Corollary \ref{cor:mixed}. \dia 
\end{ex} 

\noindent 
In our preceding $3\times 3$ example, each parallelepiped corresponds 
naturally to a $3\times 3$ {\em binomial} system, easily 
obtainable from each triplet of edges mentioned above. So the 
intersections of the $\atrop(f_i)$ naturally yield approximations to 
complex {\em roots} of $F$ (not just their norms). This approach to 
canonical start points for homotopy continuation is pursued further in 
\cite{aggr}. See \cite{phrojas} for the relevance of the preceding system 
to fewnomial theory over general local fields. 

\subsection{On the Computational Complexity of $\pmb{\atrop(f)}$ and 
$\pmb{\amoeba(f)}$}  
\label{sub:cxity} 
The\linebreak 
complexity classes $\pp$, $\np$, $\pspa$, and $\expt$ --- from the 
classical Turing model of computation --- can be identified 
with families of {\em decision} problems, i.e., problems with a yes or no 
answer. Larger complexity classes correspond to problems with larger 
worst-case complexity. We refer the reader to \cite{sipser,papa,arora,sipser2} 
for further background. Aside from the basic definitions of input size and 
$\np$-hardness, it will suffice here to simply recall that 
$\pp\!\subseteq\!\np\!\subseteq\!\pspa\!\subseteq\!\expt$, and that the 
properness of each inclusion (aside from $\pp\!\subsetneqq\!\expt$, which 
has been known for some time \cite{hartmanis,robson}) is a famous open problem. 
All algorithmic complexity results below count bit operations, and do so as a 
function of some underlying notion of input size. 

Deciding membership in an amoeba can easily be rephrased as a problem 
within the\linebreak Existential Theory of the Reals. The latter setting 
has been studied extensively in the 20$\thth$ century (see, e.g., 
\cite{tarski,bkr,canny}) and the current state 
of the art implies that amoeba membership is in $\pspa$, i.e., it 
can be solved in polynomial-time by a parallel algorithm\footnote{i.e., an 
algorithm distributed across several 
processors running simultaneously on some shared memory...}, 
provided one allows exponentially many processors. More precisely, we define 
the {\em input size} of a polynomial $f\!\in\!\Z[x_1,\ldots,x_n]$, written 
$f(x)\!=\!\sum^t_{i=1} 
c_i x^{a_i}$, to be 
$\size(f)\!:=\!\sum^t_{i=1} \log_2\left((2+|c_i|)\prod^n_{j=1}(2+|a_{i,j}|)
\right)$, where $a_i\!=\!(a_{i,1},\ldots,a_{i,n})$ for all $i$. 
(Put another way, up to a constant additive error, $\size(f)$ is just the sum 
of the bit-sizes of all the coefficients and exponents.) Similarly, we define 
$\size(v)$, for any $v\!=\!(v_1,\ldots,v_n)\!\in\!\Qn$, to be the 
sum of the sizes of the numerators and denominators 
of the $v_i$ (written in lowest terms). We similarly extend the notion of 
input size to polynomials in $\Q[x_1,\ldots,x_n]$. Considering real and 
imaginary parts, we can extend further still to 
polynomials in $\Q\!\left[\sqrt{-1}\right]\![x_1,\ldots,x_n]$. 
\begin{rem} 
For our notion of input size, sufficiently sparse polynomials have 
size\linebreak polynomial in the {\em logarithm} of the degree of the
polynomial (among other parameters). For instance, our definition implies that 
$c+x_1+x^d_2$ has size $O(\log(c)+\log(d))$. This is in contrast to 
other definitions of input size in older papers (see, e.g., \cite{theobald}) 
where degree is counted in such a way that $1+x_1+x^d_2$ has size $\geq\!d$. 
\dia 
\end{rem} 
\begin{thm} 
\label{thm:cxity} 
There is a $\pspa$ algorithm to decide, for any input 
pair $(z,f)\in$\linebreak $\bigcup_{n\in\N}\left(\Qn\times 
\Q\!\left[\sqrt{-1}\right]\![x_1,\ldots,x_n]\right)$, 
whether $\Log|z|\!\in\!\amoeba(f)$.  Furthermore, the special 
case where $z\!=\!1$ and $f\!\in\!\Z[x_1]$ 
in the preceding membership problem is already $\np$-hard.  
\end{thm} 

\noindent 
Theorem \ref{thm:cxity} is implicit in the papers 
\cite{plaisted,bkr,canny} so, for the convenience of the reader, we 
provide an outline of the proof in Section \ref{sub:baxity}. 
\begin{rem} 
While \cite[Cor.\ 2.7]{theobald} mentions ``polynomial-time'' 
amoeba membership\linebreak detection in fixed dimension, the definition of 
input size implicitly used in \cite{theobald} differs from ours and yields 
complexity polynomial in the degree, among other parameters. So the method 
underlying \cite[Cor.\ 2.7]{theobald} in fact has {\em exponential} worst-case 
complexity relative to the input size we use here. Indeed, 
the $\np$-hardness lower bound from Theorem \ref{thm:cxity} 
tells us that the existence of a polynomial-time amoeba membership algorithm 
for $n\!=\!1$ (relative to our notion of input size here) would imply 
$\pp\!=\!\np$. \dia 
\end{rem} 

Since we now know that $\atrop(f)$ is provably close to $\amoeba(f)$, 
$\atrop(f)$ would be of great practical value if $\atrop(f)$ were easier to 
work with than $\amoeba(f)$. This indeed appears to be the case. For example, 
when the dimension $n$ is fixed and all the coefficient absolute values of 
$f$ have rational logarithms, standard high-dimensional convex hull algorithms 
(see, e.g., \cite{edelsbrunner}) enable us to  
describe every face of $\atrop(f)$, as an explicit 
intersection of half-spaces, in polynomial-time (see, e.g., \cite{aggr}).  

\scalebox{.94}[1]{The case of rational coefficients presents some subtleties 
because the underlying computations,}\linebreak 
\scalebox{.915}[1]{done naively, involve arithmetic on rational numbers 
with exponentially large bit-size. Nevertheless,}\linebreak 
point membership for $\atrop(f)$ has polynomial bit complexity when $n$ is 
fixed. 
\begin{thm} 
\label{thm:cxity2} 
Fix any $\eps\!>\!0$. Then there is an 
$O\!\left(nt(\log d)^{1+\eps} (20.8\sigma(\log \sigma)^{1+\eps})^{2n+2}
\right)$\linebreak  
algorithm to decide, for any 
input $(z,f)\!\in\!\bigcup_{n\in\N}\left(\Q^n\times 
\Q\!\left[\sqrt{-1}\right]\![x_1,\ldots,x_n]\right)$ (with 
$f(x)$ of the form $\sum^t_{i=1} c_i x^{a_i}$, with degree at most $d$ with 
respect to any variable, $z\!=\!(z_1,\ldots,z_n)$, and the bit-sizes
of the $z_i$ and $c_i$ at most $\sigma$), whether $\Log|z|\!\in\!\atrop(f)$. 

Furthermore, if we instead assume that both $\log|z_i|,\log|c_i|\!\in\Q$ 
have bit size $\leq\!\sigma$ for all $i$, then there is an 
$O\!\left(nt(\sigma+\log d)\log^2(\sigma d) \right)$ algorithm to decide 
whether $\Log|z|\!\in\!\atrop(f)$.  
\end{thm} 

\noindent 
We prove Theorem \ref{thm:cxity2} in Section \ref{sec:cxity2}. 
The complexity of finding the distance to $\atrop(f)$ from a given query 
point $v$, and the relevance of such distance computations  
to polynomial system solving, is explored further in \cite{aggr}. 

\subsection{Non-Archimedean Precursors and Simplified Maslov Dequantization} 
\label{sub:arch} \mbox{}\linebreak 
Recall that $\C\langle\langle s\rangle \rangle$ is the 
union of formal Laurent series fields
$\bigcup_{d\in \N} \C\left(\left(s^{1/d}\right)\right)$. 
While $\C$ is\linebreak 
\scalebox{.98}[1]{perhaps a more popular field in applications than $\C\langle
\langle s\rangle\rangle$, $\C$ is more exceptional algebraically:}\linebreak 
$\C$ is the unique (up to isomorphism) algebraically closed field that is 
complete with\linebreak 
respect to an absolute value that is unbounded on $\Z$ (see, e.g., 
\cite[Thm.\ 1.2.3]{valued}). Such an\linebreak 
\scalebox{.95}[1]{absolute value is called 
{\em Archimedean}, so let us now review what a {\em non}-Archimedean 
valuation is.}  

A {\em (non-Archimedean) valuation} on a field $K$ is  
a function $\nu : K \longrightarrow 
\R\cup\{\infty\}$ such that $\nu(0)\!=\!\infty$ and, for all 
$a,b\!\in\!K$, (1) $\nu(ab)\!=\!\nu(a)+\nu(b)$ and (2) 
$\nu(a+b)\!\geq\!\min\{\nu(a),\nu(b)\}$ with equality if 
$\nu(a)\!\neq\!\nu(b)$. Inequality (2) is sometimes called the 
{\em Ultrametric Inequality}. Note in particular that 
$-\log(a+b)\!\geq\!\min\{-\log a,-\log b\}-\log 2$ for any 
$a,b\!\in\!\R_+$, so $-\log|\cdot|$ violates the Ultrametric Inequality 
when $a\!=\!b\!=\!1$. $-\log|\cdot|$ is thus sometimes called\linebreak   
\scalebox{.99}[1]{the {\em Archimedean} valuation on $\C$, and the minus sign 
on the $\log$ is one of the reasons behind}\linebreak 
\scalebox{.94}[1]{various sign discrepancies when comparing 
Archimedean and non-Archimedean tropical varieties.}  

More to the point, let us recall a particular classical non-Archimedean 
valuation used often in the current tropical geometry literature. 
\begin{dfn} 
\scalebox{.95}[1]{We define the {\em $s$-adic valuation}
of any element $c\!=\!\sum^\infty_{j=k} \gamma_j s^{j/d}\in 
\C\langle\langle s\rangle \rangle\setminus\{0\}$}\linebreak to be 
$\nu_s(c)\!:=\!\min_{\gamma_j\neq 0} j/d$, and set $\nu_s(0)\!:=\!\infty$.  
We then define the {\em $s$-adic Newton polytope} of any 
$f\!\in\!\C\langle \langle s\rangle \rangle\!\left[x^{\pm 1}_1,
\ldots,x^{\pm 1}_n\right]$ to be 
$\newt_s(f)\!:=\!\conv\!\left(\{(a_i,\nu_s(c_i))\}_{i\in [t]}\right)$.
We also define the {\em $s$-adic tropical variety of $f$} to be\\ 
\scalebox{.97}[1]{$\trop_s(f)\!:=\!\{v\!\in\!\Rn\; | \; 
(v,1) \text{ is an inner normal of a face of } \newt_s(f) \text{ of 
positive dimension}\}$. \dia}  
\end{dfn} 
\begin{ex}
We have drawn $\newt_s(f)$ for the trinomial $f(x_1)\!:=\!s-x^{16}_1+x^{49}_1$ 
below, along with some representative inner normals for the 
edges of $\newt_s(f)$:\\ 
\scalebox{.95}[.95]{\begin{picture}(500,30)(0,-10)
\put(0,10){{\color{black}\circle*{3}}}
\put(0,10){\line(16,-1){160}}
\put(160,0){{\color{black}\circle*{3}}}
\put(160,0){\line(1,0){330}}
\put(490,0){{\color{black}\circle*{3}}}
\put(490,0){\line(-49,1){490}}
\put(70,-7){\vector(1,16){.75}}
\put(320,-13){\vector(0,1){12}}
\put(240,18){\vector(-1,-49){.2448}}
\end{picture}}
There are just two upward-pointing inner normals, and thus just
two inner normals of the form $(v,1)$: $\left(\frac{1}{16},1\right)$ and 
$(0,1)$. So $\trop_s(f)\!=\!\left\{\frac{1}{16},0\right\}$ here. \dia 
\end{ex} 

Letting $Z^*_s(f)$ denote the roots of $f$ in $(\C\langle\langle s\rangle 
\rangle\setminus\{0\})^n$, Newton's 17$\thth$-century result on 
Puiseux series expansions \cite{newton}, in modern language, can then be 
paraphrased as follows: 
\begin{thm} 
\label{thm:newton} 
\cite{newton} 
If $f\!\in\!\C\langle \langle s\rangle \rangle\!\left[x^{\pm 1}_1
\right]$, $v\!\in\!\Q$, and $L$ is the face of $\newt_s(f)$ with 
inner normal $(v,-1)$, then $f$ has exactly $\lambda(L)$ roots, 
counting multiplicity, with $s$-adic valuation $v$. In particular, 
$\nu_s\!\left(Z^*_s(f)\right)\!=\!\trop_s(f)$. \qed 
\end{thm} 
\begin{ex} 
The trinomial $f(x_1)\!:=\!s-x^{16}_1+x^{49}_1$ from 
our last example has exactly
$49$ roots in $\C\langle\langle s\rangle \rangle$: $16$ of the 
form $e^{2\pi \sqrt{-1}j/16}s^{1/16}+\sum^\infty_{i=2}\alpha_{i,j}s^{i/16}$
(for $j\!\in\![16]$) and $33$ of the
form $e^{2\pi \sqrt{-1} j/33}+\sum^\infty_{i=1}\beta_{i,j}s^i$ 
(for $j\!\in\!
[33]$), where $\alpha_{i,j}\!\in\!\Q\!\left[e^{2\pi \sqrt{-1}/16}\right]$
and $\beta_{i,j}\!\in\!\Q\!\left[e^{2\pi \sqrt{-1}/33}\right]$. 
So the horizontal lengths ($16$ and $33$) of the two lower edges of 
$\newt_s(f)$ indeed count the number of roots with corresponding 
valuation. \dia
\end{ex}

Note how the valuations $\nu_s\!\left(Z^*_s(f)\right)$ are 
{\em exactly} determined by the lower 
edges of $\newt_s(f)$, unlike the Archimedean setting where approximation 
is unavoidable (witness Example \ref{ex:tri}). Our Theorem \ref{thm:harder} 
is thus an Archimedean analogue of Newton's result, including 
counting norms up to some variant
of multiplicity. Dumas, around 1906, extended Theorem 
\ref{thm:newton} to the $p$-adic complex numbers $\C_p$ \cite{dumas}. In 
fact, one can replace $\C\langle \langle s\rangle\rangle$ by any 
algebraically closed field with non-Archimedean valuation \cite{weiss}. 

There are two important additional characterizations of 
$\trop_s(f)$ (resp.\ $\atrop(f)$). 
\begin{prop}
\label{prop:slope}
For any univariate $f\!\in\!\C\langle\langle s\rangle\rangle\!\left[
x^{\pm 1}_1\right]$ (resp.\ 
$f\!\in\!\C\!\left[x^{\pm 1}_1\right]$) we have that $-\trop_s(f)$ (resp.\ 
$\atrop(f)$) is the set of slopes of the lower edges of $\newt_s(f)$ (resp.\ 
$\anewt(f)$). \qed 
\end{prop}
\begin{lemma}
\label{lemma:trop} 
For any $f\!\in\!\C\langle \langle s\rangle \rangle\!\left[x^{\pm 1}_1,
\ldots,x^{\pm 1}_n\right]$ we have that\\  
\mbox{}\hfill $
\displaystyle{\trop_s(f)\!=\!\left\{v\!\in\!\Rn\; \left| \; 
\min_i\{a_i\cdot v+ \nu_s(c_i)\} \text{ is attained for at least two distinct 
values of } i\right.\right\}}$.\hfill\mbox{}\\ 
Also, for any $f\!\in\!\C\!\left[x^{\pm 1}_1,
\ldots,x^{\pm 1}_n\right]$ we have that\\
\mbox{}\hfill $
\displaystyle{\atrop(f)\!=\!\left\{v\!\in\!\Rn\; \left| \; 
\max_i|c_ie^{a_i\cdot v}| \text{ is attained for at least two distinct values
of } i\right.\right\}}$.\hfill\qed 
\end{lemma} 

\noindent 
Proposition \ref{prop:slope} is elementary, while Lemma 
\ref{lemma:trop} follows immediately from the definition of an inner/outer 
face normal (see, e.g., \cite[Ch.\ 7]{ziegler}). In particular, in the 
non-Archimedean 
case, these alternative characterizations are well-known in the tropical 
literature (see, e.g., \cite{macsturmf}). In the Archimedean case, 
Hadamard and Ostrowski's original univariate root norm estimates were in fact 
stated in terms of edge slopes.  

For any $\zeta_1,\ldots,\zeta_n\!\in\!\C\langle\langle s\rangle 
\rangle\setminus\{0\}$, let $\nu_s(\zeta)\!:=\!(\nu_s(\zeta_1),
\ldots,\nu_s(\zeta_n))$. That the last assertion of Theorem \ref{thm:newton} 
can be extended to multivariate polynomials was first observed by Kapranov. 
\begin{kapra} (Special Case)\footnote{Kapranov's Theorem was originally 
stated for any algebraically closed field with a non-Archimedean 
valuation. Note in particular that $\nu_s(\C\langle\langle s
\rangle\rangle)\!=\!\Q$ here. } \cite{ekl} 
For any $f\in$\linebreak 
$\C\langle\langle s\rangle\rangle\!
\left[x^{\pm 1}_1,\ldots,x^{\pm 1}_n\right]$,   
we have $\nu_s\!\left(Z^*_s(f)\right)=\trop_s(f)\cap \nu_s(\C\langle\langle 
s \rangle\rangle)^n$. \qed 
\end{kapra}

In Kapranov's Theorem (as well as Theorem \ref{thm:newton}), the containment 
of $s$-adic root valuation vectors in $\trop_s(f)\cap \nu_s(\C\langle\langle 
s \rangle\rangle)^n$ follows easily from the Ultrametric Inequality. 
Over the Archimedean field $\C$, proving that $\amoeba(f)$ is contained in a 
suitable neighborhood of $\atrop(f)$ requires a more delicate application of 
the  Triangle Inequality. Proving that $\atrop(f)$ is contained in a suitable 
neighborhood of $\amoeba(f)$ then involves specializing\linebreak 
\scalebox{.965}[1]{to a curve (similar to 
a trick in the non-Archimedean setting) to reduce to the univariate case,}
\linebreak   
and then applying Rouch\'e's Theorem. An estimate on lattice points 
visible from the origin (Theorem \ref{thm:whew} in Section \ref{sub:big} 
below) helps improve one of our bounds in the bivariate case. 

We close with some topological observations. First observe that $\atrop(f)$ 
need {\em not} be contained in $\amoeba(f)$, nor even have the same homotopy 
type as $\amoeba(f)$,\linebreak already for $n\!=\!1$: The example 
$f(x_1)\!=\!(x_1+1)^2$ yields $\atrop(f)\!=\!\{\pm \log 2\}$ but 
$\amoeba(f)\!=\!\{0\}$. However, one can in fact always recover $\atrop(f)$ 
as the Hausdorff limit of a {\em sequence of suitably scaled} amoebae. To 
clarify this, first recall that the {\em Hausdorff}\linebreak {\em distance} 
between any two subsets $X,Y\!\subseteq\!\Rn$ is \\ 
\mbox{}\hfill 
$\Delta(X,Y)\!:=\!\max\left\{\sup\limits_{x\in X}{} 
\inf\limits_{\substack{\mbox{}\\ y\in Y}}|x-y|,
\sup\limits_{y\in Y}{} \inf\limits_{\substack{\mbox{}\\ x\in X}}
|x-y|\right\}$.\hfill\mbox{}\\ 
Also, the {\em support} of a Laurent polynomial 
$f(x)\!=\!\sum^t_{i=1} c_ix^{a_i}$ is 
$\supp(f)\!:=\!\{a_i\; | \; c_i\!\neq\!0\}$. 
\begin{cor} 
\label{cor:hausdorff} Let $f$ be any $n$-variate 
$t$-nomial with $t\!\geq\!2$ and $k\!:=\!\dim\newt(f)$. Then:  
\begin{enumerate} 
\item{$\Delta(\atrop(f),\amoeba(f))\!\leq\!\sqrt{k}\ceil{\frac{1}{4} 
t(t-1)}\left((\log 9)t-\log\frac{81}{2}\right)\!=\!O\!\left(
t^{7/2}\right)$.} 
\item{There exists a family of Laurent polynomials $(f_\mu)_{\mu\geq 1}$ 
with $\supp(f_\mu)\!=\!\supp(f)$ for all $\mu\!\geq\!1$ and 
$\Delta\!\left(\frac{1}{\mu}\amoeba(f_\mu),\atrop(f)\right)
\longrightarrow 0$ as $\mu\longrightarrow \infty$.}  
\end{enumerate} 
\end{cor} 

\noindent  
We will prove Corollary \ref{cor:hausdorff} momentarily, but 
let us first recall one of the consequences of {\em Maslov dequantization} 
(see, e.g., \cite{maslov,lms,virologpaper} and \cite[Cor.\ 6.4]{mik3}): A 
method to obtain any {\em non}-Archimedean tropical variety as a limit of a 
family of scaled Archimedean amoebae. Assertion (2) thus shows how $\atrop(f)$ 
provides a fully Archimedean version of this limit. 
Another precursor to Assertion (2), involving the piecewise linear 
structure approached by the intersection of $\amoeba(f)$ with a large 
sphere, appears in \cite{bergman} and \cite[Prop.\ 1.9, Pg.\ 197]{gkz94}. 
Thanks to Assertion (1), we can prove Assertion (2) in just three lines. 

\medskip
\noindent 
{\bf Proof of Corollary \ref{cor:hausdorff}:} 
Assertion (1) of Corollary \ref{cor:hausdorff} follows immediately 
from\linebreak 
Assertion (1b) of Theorem \ref{thm:finally}, and the fact that 
$k\!\leq\!t-1$. Let us write 
$f(x)\!=\!\sum^t_{i=1} c_ix^{a_i}$, define $f_\mu(x)
\!:=\!\sum^t_{i=1} c^\mu_i x^{a_i}$, and observe that 
$f_1\!=\!f$.  

Since $|c_ie^{a_i\cdot v}|\!\geq\!|c_je^{a_j\cdot v}| 
\Longleftrightarrow 
|c_ie^{a_i\cdot v}|^\mu\!\geq\!|c_je^{a_j\cdot v}|^\mu$,
we immediately obtain that $\atrop(f_\mu)$\linebreak 
\scalebox{.935}[1]{$=\!\mu\atrop(f)$. So then
$\Delta\!\left(\amoeba(f_\mu),\atrop(f_\mu)\right) 
=\!\mu\Delta\!\left(\frac{1}{\mu}\amoeba(f_\mu),\atrop(f)
\right)$ and}\linebreak Assertion (1) thus implies 
$\Delta\!\left(\frac{1}{\mu}\amoeba(f_\mu),\atrop(f)\right)
\!=\!\frac{O\left(t^{7/2}\right)}{\mu}$ for all $\mu\!\geq\!1$. \qed 

\section{Background on Univariate Bounds and the Complexity of 
Amoeba Membership} 
\label{sec:back} 
To prepare for the proofs of our main metric results we will first review some 
classical root norm bounds in the univariate case, in order to recast them in 
terms of $\atrop(f)$. We then prove 
Theorem \ref{thm:harder} (along with a refinement for $t\!=\!3$), 
and conclude this section with a sketch of the proof of Theorem 
\ref{thm:cxity} (on the hardness of deciding point membership for amoebae). 
\begin{rem} 
\label{rem:nota} 
Throughout this section we assume that $f(x_1)\!=\!\sum^t_{i=1} c_i 
x^{a_i}\!\in\!\C\!\left[x^{\pm 1}_1\right]$ is a univariate $t$-nomial with 
$t\!\geq\!2$ (so the $c_i$ are all nonzero), $0\!\leq\!a_1\!<\cdots<\!a_t$, 
and $d\!:=\!a_t-a_1$. We then label the roots of $f$ in $\Cs$ by  
$\zeta_1,\ldots,\zeta_d$, counting multiplicity, in such a way that 
$|\zeta_1|\!\leq\cdots\leq\!|\zeta_d|$. For each $i\!\in\!\{1,\ldots,d\}$ 
we also let $v_i$ denote the slope of the (unique) lower edge of the polygon  
$\anewt(f)\cap ([a_1+i-1,a_1+i]\times \R)$. Since each such lower edge 
is in fact a line segment inside a lower edge of $\anewt(f)$, 
Proposition \ref{prop:slope} thus implies that 
$\{v_1,\ldots,v_d\}\!=\!\atrop(f)$. We also 
order the $v_i$ so that $v_1\!\leq\cdots \leq\!v_d$. \dia 
\end{rem} 

We begin with a pair of bounds dating back to 1893 and 1923. 
\begin{thm}
\label{thm:old} Following the notation of Remark \ref{rem:nota}, the following 
hold:\\ 
(1) If $\zeta\!\in\!\C$ is a root of $\alpha_0+\cdots+\alpha_d x^d_1\!\in\!
\C\!\left[x^{\pm 1}_1\right]$, and $\alpha_0\alpha_d\!\neq\!0$, then\\ 
\mbox{}\hfill  
$\frac{1}{2}\min\limits_{\alpha_i\neq 0,i\neq 0}\left|\frac{\alpha_0}{\alpha_i}
\right|^{1/i} < |\zeta| < 2\max\limits_{i\in\{0,\ldots,d-1\}}
\left|\frac{\alpha_i}{\alpha_d}\right|^{1/(d-i)}$.\hfill\mbox{}\\  
(2) If $g(x_1)\!=\!\beta_0+\cdots+\beta_px^p_1+ 
\gamma_1x^{n_1}_1+\cdots+\gamma_q x^{n_q}_1\!\in\!\C[x_1]$ so that   
$\beta_p\!\neq\!0$ and\\  
\mbox{}\hspace{.6cm}$1\!\leq\!p\!<\!n_1\!<\cdots<\!n_q$, 
then $g$ has a root with absolute value 
$\leq\!\left|\frac{\beta_0}{\beta_p}\right|^{1/p} \binom{p+q}{q}^{1/p}$. \qed 
\end{thm} 

\noindent 
Bound (1) is a paraphrase of a special case of \cite[Pg.\ 
201, Third Inequality]{hadamard}, and is stated more explicitly in 
\cite{fujiwara}   
(see also \cite[pp.\ 243--249]{rs}, particularly Bound 8.1.11 on Pg.\ 247).   
Bound (2) was proved by Montel \cite{montel} (see also \cite[Thm.\ 9.5.1, 
Pg.\ 304]{rs}). 
\begin{cor} 
\label{cor:extremecauchy} Following the notation of Remark \ref{rem:nota}, 
we have:\\ 
\mbox{}\hspace{2.1cm}(a)\hspace{1.05cm}$-\log 2< \log|\zeta_1|-\min\atrop(f) \, 
\leq \log(t-1)$\\ 
\mbox{}\hspace{2.1cm}(b) $-\log(t-1)\leq \log|\zeta_d|-\max\atrop(f) < 
\log 2$\\
\mbox{}\hspace{2.1cm}(c) The $\log 2$ (resp.\ $\log(t-1)$) terms above can not 
be replaced by any smaller\\ 
\mbox{}\hspace{2.8cm}constant (resp.\ function of $t$ solely).  
\end{cor} 

\noindent 
{\bf Proof:} The lower bound from Part (a) and the upper bound 
from Part (b) follow immediately from Proposition \ref{prop:slope},  
upon taking the log absolute value of both sides of Bound (1) from Theorem 
\ref{thm:old}. In particular, we see that the lower and upper bounds 
from Bound (1) are exactly $\frac{1}{2}e^{\min \atrop(f)}$ and 
$2e^{\max \atrop(f)}$. 

The upper bound from Part (a) follows similarly, but employing Bound (2) from 
Theorem \ref{thm:old} instead of Bound (1). In particular, one 
must apply Bound (2) in the following way: Take $p$ so that 
the $(p,-\log|\beta_p|)$ is the right-hand vertex of the left-most 
lower edge of $\anewt(f)$. By construction, this edge has slope 
$\frac{\log|\beta_0|-\log|\beta_p|}{p}$. Observing that 
$\binom{p+q}{q}^{1/p}=$\linebreak 
$\left(\frac{(q+p)\cdots (q+1)}{p!}\right)^{1/p}\!=\!
\left(\left(\frac{q}{p}+1\right)\cdots \left(\frac{q}{1}+1\right)
\right)^{1/p}\!\leq\!\left((q+1)^p\right)^{1/p}\!=\!q+1$, and that the 
number of terms is $t\!=\!p+q+1$ with $p\!\geq\!1$, we are done. 

The lower bound from Part (b) follows by applying the preceding 
paragraph to the polynomial $x^{a_t+a_1}_1f(1/x_1)$: This has the effect 
of reflecting $\anewt(f)$ across the vertical line $\frac{d}{2}\times \R$, 
and thus $\atrop(f)$ is replaced by $-\atrop(f)$. So we ultimately 
prove an upper bound of $\log(t-1)$ on $-\log|\zeta_d|-(-\max\atrop(f))$ 
and we are done. 

The optimality of the $\log 2$ terms is evinced by the polynomials\\  
\mbox{}\hfill 
$f_1(x_1)\!:=\!x^{t-1}_1-x^{t-2}_1-\cdots-1$ \ \ \ and \ \ \ 
$f_2(x_1)\!:=\!-1+x_1+ \cdots+x^{t-1}_1$:\hfill\mbox{}\\ 
One need only show that $f_1$ (resp.\ $f_2$) 
has a unique positive root increasing toward a limit of $2$ (resp.\ decreasing 
toward a limit of $\frac{1}{2}$) as 
$t\longrightarrow \infty$. Uniqueness follows from Descartes' Rule, 
since each $f_1$ and $f_2$ have exactly one sign alternation in their 
ordered sequence of coefficients. The limiting behavior of their unique 
positive roots is easily obtained from Rolle's Theorem (since 
$f_1(0),f_2(0)\!<\!0$ and $f_1(x_1),f_2(x_1)\longrightarrow \infty$ 
as $x_1\longrightarrow\infty$), and the fact that $1\!=\!\frac{1}{2}
+\frac{1}{4}+\frac{1}{8}+\cdots$. 

The optimality of the $\log(t-1)$ terms 
is easily seen via the polynomial\linebreak 
$g(x_1)\!:=\!(x_1+1)^{t-1}$: The left-most (resp.\ right-most) 
lower edge of $\anewt(g)$ has slope $-\log(t-1)$ (resp.\ 
$\log(t-1)$), by the log-concavity of the binomial coefficients. 
So by Proposition \ref{prop:slope}, $\min \atrop(g)\!=\!
-\log(t-1)$ and $\max \atrop(g)\!=\!\log(t-1)$. Since $\amoeba(g)\!=\!\{0\}$, 
we are done. \qed 

\medskip 
By a slight variation of our last proof, we can easily obtain a family of 
examples showing that the $\log 3$ interval-width from Theorem \ref{thm:harder} 
is in fact optimal. 
\begin{lemma} 
\label{lemma:opt} For any $k\!\geq\!1$ 
let $g_k(x_1)\!:=\!1+x+\cdots+x^{k-1}-x^k+\frac{1}{9}x^{k+1} 
+\cdots+\frac{1}{9^k}x^{k+k}$. Then (1) $\atrop(g_k)\!=\!\{0,\log 9\}$ for 
all $k$ and (2) for any fixed $\eps\!>\!0$ 
there is a $k$ such that $\amoeba(f)\cap(\log(3)-\eps,\log(3)+\eps)$ is 
non-empty. \qed 
\end{lemma} 

We now recall a seminal collection of bounds due to Ostrowski: 
\begin{thm} 
\label{thm:ostrowski} 
\cite[Cor.\ IX, Pg.\ 143]{ostrowski}\footnote{There was a typo 
in Ostrowski's original statement of the upper bound from 
Assertion (3), later corrected in an addendum by Ostrowski \cite{ostrowski2}.}  
Following the notation of Remark \ref{rem:nota}, we have:\\  
\mbox{}\hspace{1.6cm}(1) \hspace{1.3cm} 
 $-\log 2< \log|\zeta_1|-v_1 \leq \log d$,\\
\mbox{}\hspace{1.6cm}(2) \hspace{1.3cm} 
 $-\log d\leq \log|\zeta_d|-v_d< \log 2$,\\
\mbox{}\hspace{1.6cm}(3) $\displaystyle{\log\left(1-\frac{1}{2^{1/i}}\right)< 
\log|\zeta_i|-v_i<-\log\left(1-\frac{1}{2^{1/(d-i+1)}}\right)}$ for 
all $i\!\in\!\{2,\ldots,d-1\}$.\\ 
In particular, \ $\displaystyle{-0.5348 
\leq \ \ \ \ \ \ \ \log\left(1-\frac{1}{2^{1/i}}\right) \; \; \; \, - \, \ \ 
(-\log i) \ \ \ \ \; < -0.3665}$ and\\ 
\mbox{}\hspace{3cm}$\displaystyle{ \, 0.3665<-\log\left(1-\frac{1}
{2^{1/(d-i+1)}}\right)-\log(d-i+1) \; \leq \ \; 0.5348}$. \hfill\qed\\ 
\end{thm} 
\begin{rem}
Since $\atrop(f)\!=\!\{v_1,\ldots,v_d\}$ (see Remark \ref{rem:nota}), 
our Theorem \ref{thm:harder} implies that any given 
$\log|\zeta_i|$ lies within distance $\log 3$ of some $v_j$, {\em possibly 
with $j\!\neq\!i$}. 
In this sense, the final assertion of Theorem \ref{thm:ostrowski} 
tells us that Theorem \ref{thm:harder} isolates each $\log|\zeta_i|$ strictly 
better than Ostrowski's bounds, except possibly in the cases  
$i\!\in\!\{2,d-1\}$ or $t\!=\!d+1\!=\!3$. Corollary \ref{cor:tighter} 
in Section \ref{sub:harder} below matches Ostrowski's bounds when 
$t\!=\!d+1\!=\!3$. \dia 
\end{rem} 

More recently, Akian, Gaubert, and Sharify have derived metric 
bounds improving those of Hadamard and Ostrowski \cite{ags}, but in a different 
direction from ours. Their focus was the {\em Matrix Polynomial 
Problem} (a.k.a.\ the {\em Polynomial Eigenvalue Problem}): 
Given matrices $A_0,\ldots,A_d\!\in\!\C^{n\times n}$, find 
$\lambda\!\in\!\C$ such that $A_0+\lambda A_1+\cdots+\lambda^d A_d$ 
has determinant $0$. The Matrix Polynomial Problem includes the 
classical eigenvalue problem for $d\!=\!1$, while the $n\!=\!1$ case 
is the problem of univariate polynomial solving. A special case of one 
of the main theorems of \cite{ags} adds a new bound to the univariate 
$t$-nomial setting: The 
center of mass of the $\ell$ smallest points of $\amoeba(f)$ is not 
too far to the left of the center of mass of the $k$ smallest points of 
$\atrop(f)$ --- assuming repeated points are counted appropriately. 
\begin{thm} 
\label{thm:ags} 
(Special case of \cite[Thm.\ 4.1 \& Rem.\ 4.2]{ags}) Following the 
notation of Remark \ref{rem:nota}, for any $\ell\!\in\!\{1,\ldots,d\}$, 
we have $\sum^\ell_{i=1}\log|\zeta_i|\!\geq\!\left(
\sum^\ell_{i=1}v_i\right)-\frac{1}{2}\log t$. \qed 
\end{thm} 
\begin{rem} 
\scalebox{.93}[1]{It is worth observing some notational divergences: In 
the notation of \cite{ags}, our}\linebreak 
\scalebox{.87}[1]{$\atrop(f)$ would be the log of the set of ``tropical roots'' 
of the ``tropical polynomial'' $\max_{i\in\{1,\ldots,t\}} 
|c_i|x^{a_i}$}\linebreak 
defined over the nonnegative reals. However, most other authors (e.g., 
\cite{pin,macsturmf}) would instead call 
$\max_{i\in \{1,\ldots,t\}} \{a_i w + \log|c_i|\}$ or 
$\min_{i\in \{1,\ldots,t\}} \{a_i w - \log|c_i|\}$ 
a tropical polynomial, depending on what semi-ring they prefer. 
This in turn introduces a sign flip in their corresponding 
definition of tropical root or tropical variety. 
(The oldest definition of tropical polynomial is in fact via 
the minimum of a collection of linear forms \cite{pin}.) 
Because of this sign discrepancy, some authors (\cite{ags} included) use a 
variant of $\anewt(f)$  which implies examining {\em upper} hulls, instead of 
lower hulls, in order to define tropical varieties. \dia  
\end{rem} 

\noindent 
Note that while the $\ell\!=\!1$ case of Theorem \ref{thm:ags} 
yields a weaker bound than Assertion (1) of Theorem \ref{thm:ostrowski} 
when $t\!\geq\!4$, a strength of Theorem \ref{thm:ags} is its new bound on a 
particular 
{\em amortized} error of approximating $\amoeba(f)$ by $\atrop(f)$: Applying 
Theorem \ref{thm:ags} with $\ell\!=\!d$ to both $f(x_1)$ and $f(1/x_1)$ implies 
that $\left|\frac{\sum^d_{i=1}\log|\zeta_i|}{d}-\frac{\sum^d_{i=1}v_i}{d}
\right|\!\leq\!\frac{1}{2d}\log t$. Such an estimate does not appear to be 
directly obtainable from our methods here. 
\begin{rem} 
\scalebox{.95}[1]{The Matrix Polynomial 
Problem can be naturally phrased as a polynomial system}\linebreak   
with solutions in $\C\times \Pro^{n-1}_\C$ 
(i.e., an {\em intersection of several} hypersurfaces)  
by considering the vector 
equality $(A_0+\lambda A_1+\cdots+\lambda^d A_d)x\!=\!\bO$. 
So the metric bounds of \cite{ags} for the case $n\!\geq\!2$ (which have 
exponential dependence on the dimension $n$) are not directly 
comparable\linebreak 
with Theorem \ref{thm:finally} (which applies to a {\em single} hypersurface, 
and has dependence sub-cubic in the number of terms). \dia  
\end{rem} 

We have so far concentrated on showing that each $\log|\zeta_i|$ 
is close to some $v_j$, with optimal distance bounds. Showing that each 
$v_j$ is close to some $\log|\zeta_i|$ requires more\linebreak 
preparation, which we now detail. 

\subsection{Proving Theorem \ref{thm:harder}} 
\label{sub:harder} 
We will need three technical results on bounding the norms of 
summands of sparse polynomials, and counting roots of polynomials 
in annuli, before proving Theorem \ref{thm:harder}. 
\begin{prop}
\label{prop:lop}
Suppose $t\!\geq\!3$, $v\!\in\!\atrop(f)$, and $\ell$ is the unique
index such that $(a_\ell,-\log|c_\ell|)$ is the right-hand vertex of
the lower edge of $\anewt(f)$ of slope $v$ (so $2\!\leq\!\ell$). Then for any 
$N\!\in\!\N$ and $x_1$ with $|x_1|\!\geq\!(N+1)e^v$ 
we have $\sum\limits^{\ell-1}_{j=1} \left|c_jx^{a_j}_1\right|\!<\!
\frac{1}{N}\left|c_\ell x^{a_\ell}_1\right|$. 
\end{prop}  

\noindent
{\bf Proof:} First note that $2\!\leq\!\ell\!\leq\!t$ by construction. 
Letting $r\!:=\!\log|x_1|$ and $\beta_j\!:=\!\log|c_j|$ we obtain 
$\sum^{\ell-1}_{j=1} \left|c_j x^{a_j}_1\right| 
   = \sum^{\ell-1}_{j=1} e^{a_jr+\beta_j}  
   = \sum^{\ell-1}_{j=1} e^{a_j(r-v)+a_jv+\beta_j}$.   
Clearly, 
$a_j\!\leq\!a_\ell-(\ell-j)$, so for $r\!\geq\!v$ we have 
$\displaystyle{
\sum^{\ell-1}_{j=1} e^{(a_\ell-(\ell-j))(r-v)+a_jv+\beta_j}\\ 
\leq \sum^{\ell-1}_{j=1} e^{(a_\ell-(\ell-j))(r-v)
                              +a_\ell v+\beta_\ell}}$, 
thanks to Proposition \ref{prop:slope} and 
the definition of $\atrop(f)$. So then\\ 
\mbox{}\hspace{1.45cm} 
$\sum^{\ell-1}_{j=1} \left|c_j x^{a_j}_1\right|\;  
\leq \; e^{(a_\ell-(\ell-1))(r-v)+a_\ell v+\beta_\ell}
        \sum^{\ell-1}_{j=1} e^{(j-1)(r-v)}$ 
\begin{eqnarray*}
& = & e^{(a_\ell-(\ell-1))(r-v)+a_\ell v+\beta_\ell}
      \left(\frac{e^{(\ell-1)(r-v)}-1}
      {e^{(r-v)}-1}\right)\\
& < & e^{(a_\ell-(\ell-1))(r-v)+a_\ell v+\beta_\ell}
      \left(\frac{e^{(\ell-1)(r-v)}}
      {e^{r-v}-1}\right) =  \frac{e^{a_\ell r+\beta_\ell}}
      {e^{r-v}-1}
\end{eqnarray*}
So to prove our desired inequality it clearly suffices to enforce
$e^{r-v}-1\!\geq\!N$. The last inequality clearly
holds for all $r\!\geq\!v+\log (N+1)$, so we are done. \qed

\medskip 
\begin{pellet} \cite{pellet}\footnote{See also 
\cite[Thm.\ 9.2.2, Pg.\ 285]{rs}. Our paraphrase here follows immediately 
from Pellet's original theorem simply by multiplying $f$ by $x^{-a_1}_1$.} 
If the Laurent polynomial 
$|c_\ell|x^{a_\ell}\;-\!\!\!\!\sum\limits_{i\in\{1,\ldots,t\}\setminus\{\ell\}} 
|c_i|x^{a_i}$ has exactly $2$ positive roots $\zeta_1\!<\!\zeta_2$ then, 
counting multiplicities, $f$ has exactly $a_\ell-a_1$ (resp.\ $a_t-a_\ell$) 
roots with norm in $(0,\zeta_1]$ (resp.\ $[\zeta_2,\infty)$). 
In particular, $f$ has no roots with norm in $(\zeta_1,\zeta_2)$. \qed  
\end{pellet} 

\begin{lemma} 
\label{lemma:vert} 
Set $\vm\!:=\!\min\atrop(f)$ and $\vp\!:=\!\max\atrop(f)$.   
Also let $v_1$ and $v_2$ be {\em consecutive} points of $\atrop(f)$ satisfying 
$v_2\!\geq\!v_1+\log 9$, and let $\ell$ be the unique
index such that $(a_\ell,-\log|c_\ell|)$ is the unique vertex of $\anewt(f)$
incident to lower edges of slopes $v_1$ and $v_2$ (so $2\!\leq\!\ell\!
\leq\!t-1$). Then, counting multiplicities, $f$ has exactly $a_\ell-a_1$ 
(resp.\ $a_t-a_\ell$) roots $\zeta\!\in\!\C$ satisfying 
$\frac{1}{2}e^{\vm}<|\zeta|<3e^{v_1}$ (resp.\ 
$\frac{1}{3}e^{v_2}<|\zeta|<2e^{\vp}$). 
\end{lemma} 

\medskip 
\noindent 
\scalebox{.95}[1]{{\bf Proof of Lemma \ref{lemma:vert}:} By symmetry (with 
respect to replacing $x_1$ by $\frac{1}{x_1}$) it clearly suffices 
to}\linebreak 
\scalebox{.95}[1]{prove the first root count. Setting $x_1\!:=\!3e^{v_1}$, 
Proposition \ref{prop:lop} tells us that $\frac{1}{2}|c_\ell| x^{a_\ell}_1\!>\!
\sum^{\ell-1}_{j=1}|c_j| x^{a_j}_1$.}\linebreak 
Observing that $\frac{1}{x_1}\!=\!\frac{1}{3e^{v_1}}\!\geq\!3e^{-v_2}$  
(since $v_2-v_1\!\geq\!\log 9$), another application of Proposition 
\ref{prop:lop} to $f(1/x_1)$ then implies that  
$\frac{1}{2}|c_\ell| x^{a_\ell}_1\!>\!
\sum^t_{j=\ell+1}|c_j| x^{a_j}_1$. 
So $g(x_1)\!:=\!|c_\ell| x^{a_\ell}_1-\sum^t_{j\neq \ell}|c_j|  
x^{a_j}_1$ is\linebreak 
positive at $x_1\!=\!3e^{v_1}$. Note also that both $g(\eps)$ and 
$g(1/\eps)$ are negative for all sufficiently small $\eps\!>\!0$. So by Rolle's 
Theorem, $g$ has at least $2$ positive roots. Moreover, by Descartes' 
Rule (since $g$ has exactly $2$ sign alternations in its ordered 
sequence of coefficients), $g$ has at most $2$ positive roots. So we may 
apply Pellet's Theorem and, applying Assertion (1) of Theorem 
\ref{thm:old} as well, we are done. \qed  

\medskip 
\noindent 
{\bf Proof of Theorem \ref{thm:harder}:} Assuming the first assertion holds, 
the final assertion is immediate from the Fundamental Theorem of Algebra: 
Counting multiplicities, $a_t-a_1$ is exactly the number of roots of $f$ 
in $\Cs$, and we clearly have 
\begin{eqnarray} 
\label{eq:1} 
\sum\limits_{\substack{\Gamma \text{ a connected component}\\ 
\text{of } \atrop(f)_{\log 3}}} \Lambda_\Gamma & = & a_t-a_1. 
\end{eqnarray} 
In particular, when $t\!\geq\!2$, 
$U_f$ is non-empty and $\Gamma$ must contain at least $1$ point of 
$\atrop(f)$. So $\Lambda_\Gamma$ is a positive integer when $t\!\geq\!2$. 

Now suppose temporarily that $\atrop(f)_{\log 3}$ is connected. Then 
$\Lambda_\Gamma\!=\!a_t-a_1$, Corollary 
\ref{cor:extremecauchy} tells us that $\amoeba(f)\!\subset\!\Gamma$, 
and we would be done. 

So assume $\atrop(f)_{\log 3}$ has at least two distinct connected components. 
Lemma \ref{lemma:vert} then immediately yields the 
conclusion of Theorem \ref{thm:harder} when $\Gamma$ is either the 
left-most or right-most connected component of $\atrop(f)_{\log 3}$: 
We simply take $v_1$ to be the right-most point of $\Gamma\cap\atrop(f)$ 
and $v_2$ the left-most point of $\atrop(f)$ in the connected component 
of $\atrop(f)_{\log 3}$ immediately to the right of $\Gamma$, or 
we take $v_2$ to be the left-most point of $\Gamma\cap\atrop(f)$ 
and $v_1$ the right-most point of $\atrop(f)$ in the connected component 
of $\atrop(f)_{\log 3}$ immediately to the left of $\Gamma$. 

\scalebox{.93}[1]{We can now proceed by 
induction on the number of connected components of $\atrop(f)_{\log 3}$:} 
\linebreak  
We simply ignore the left-most and 
right-most connected components of $\atrop(f)_{\log 3}$, and treat the 
{\em new} left-most and right-most connected components via Lemma 
\ref{lemma:vert} as in the last paragraph. 

To conclude, Corollary \ref{cor:extremecauchy} tells us that we can also 
attain $\Lambda_\Gamma$ many log norms (counting multiplicities) within the 
potentially tighter interval $\Gamma\cap (\alpha_f,\beta_f)$. 
Also, Equality (\ref{eq:1}) implies that every root of $f$ must 
have log norm within some $\Gamma$. So we are done. \qed 

\medskip 
We can tighten the union of intervals $U_f$ further when $t\!=\!3$: Combining 
Theorem \ref{thm:harder} with Assertion (1a) of Theorem \ref{thm:finally} 
(proved independently in Section \ref{sub:2}) immediately yields 
the following refinement. 
\begin{cor} 
\label{cor:tighter} 
Following the notation of Theorem \ref{thm:harder}, assume $t\!=\!3$ and 
$\beta_f-\alpha_f\!>\!\log 8$. Then there are exactly $a_2-a_1$ 
(resp.\ $a_3-a_2$) roots $\zeta\!\in\!\C$ of $f$ with 
$\alpha_f\!<\!\log|\zeta|\!\leq\!\alpha_f+\log 4$ 
(resp.\ $\beta_f-\log 4\!\leq\!\log|\zeta|\!<\!\beta_f$). \qed 
\end{cor} 

\subsection{Classical Computational Algebra and Amoeba Membership} 
\label{sub:baxity} In what follows, all $O$-constants are effective 
and absolute. Let us first recall the following results of Plaisted and 
Ben-Or, Kozen, and Reif. 
\begin{thm} 
\label{thm:plai} 
\cite{plaisted} 
The problem {\em ``Decide whether an arbitrary input $f\!\in\!\Z[x_1]$ has a 
complex root of norm $1$.''}  is $\np$-hard. \qed 
\end{thm} 
\begin{thm} 
\label{thm:canny} 
\cite{bkr,canny} 
There is an algorithm that, given any collection of polynomials
$f_1,\ldots,f_p,g_1,\ldots,g_q,h_1,\ldots,h_r\!\in\!\Q[x_1,\ldots,x_n]$, 
decides whether there is a $\zeta\!=\!(\zeta_1,\ldots,\zeta_n)\!\in\!\Rn$ 
with $f_1(\zeta)\!=\cdots=\!f_p(\zeta)\!=\!0$, $g_1(\zeta),\ldots,
g_q(\zeta)\!>\!0$, and $h_1(\zeta),\ldots,h_r(\zeta)\!\geq\!0$, 
in time\\ 
\mbox{}\hfill $\left[\sum^p_{i=1}\size(f_i))+ 
(\sum^q_{i=1}\size(g_i))+(\sum^r_{i=1}\size(h_i))\right]^{O(1)}$,\hfill 
\mbox{}\\  
using $\left[\sum^p_{i=1}\size(f_i))+ 
(\sum^q_{i=1}\size(g_i))+(\sum^r_{i=1}\size(h_i))\right]^{O(1)}$ 
processors. \qed 
\end{thm} 

\noindent 
Theorem \ref{thm:cxity} will then follow easily from two 
elementary propositions. The first is a well-\linebreak 
\scalebox{.95}[1]{known trick from 
computational algebra for re-expressing polynomial systems in a simpler form.} 
\begin{prop} 
\label{prop:1} 
Given any $f_1,\ldots,f_m\!\in\!\Q\!\left[\sqrt{-1}\right]\!
\left[x^{\pm 1}_1,\ldots,x^{\pm 1}_n\right]$, we can find 
$g_1,\ldots,g_M\!\in\!\Q\!\left[\sqrt{-1}\right]\!
\left[x^{\pm 1}_1,\ldots,x^{\pm 1}_n, y^{\pm 1}_1,\ldots,y^{\pm 1}_N
\right]$ satisfying the 
following properties: \\ 
\mbox{}\hspace{1cm}1. $f_1\!=\cdots=\!f_m\!=\!0$ 
 has a root in $\Cn \Longleftrightarrow 
g_1\!=\cdots=\!g_M\!=\!0$ has a root in $\C^N$.\\ 
\mbox{}\hspace{1cm}2. Each $g_i$ is either a quadratic binomial or a linear 
trinomial.\\ 
\mbox{}\hspace{1cm}3. $\sum^M_{i=1}\size(g_i)\!=\!O\!\left(\sum^m_{i=1}
\size(f_i)\right)$.\\ 
Moreover, $g_1,\ldots,g_M$ can be found in time 
$O\!\left(\sum^m_{i=1}\size(f_i)\right)$. \qed  
\end{prop}  
\begin{prop} 
\label{prop:2} 
Given any $f_1,\ldots,f_m\!\in\!\Q\!\left[\sqrt{-1}\right]\! 
\left[x^{\pm 1}_1,\ldots,x^{\pm 1}_n\right]$ with each $f_i$ of degree at most 
$2$, we can find $g_1,\ldots,g_M\!\in\!\Q\!\left[
x^{\pm 1}_1,\ldots,x^{\pm 1}_n, 
y^{\pm 1}_1,\ldots,y^{\pm 1}_N 
\right]$ satisfying the following properties: 1. $f_1\!=\cdots=\!f_m\!=\!0$ 
has a root in $\Cn \Longleftrightarrow 
g_1\!=\cdots=\!g_M\!=\!0$ has a root in $\R^N$.\\ 
\mbox{}\hspace{2.05cm}2. 
$\sum^M_{i=1}\size(g_i)\!=\!O\!\left(\sum^m_{i=1}\size(f_i)\right)$.\\  
Moreover, $g_1,\ldots,g_M$ can be found in time $O(mn)$. \qed 
\end{prop} 

\noindent 
A simple example of Proposition \ref{prop:1} is the replacement of 
$f(x_1)\!:=\!1-2x_1+x^5_1$ by the system $G\!:=\!(y_1-x^2_1,y_2-y^2_1,
y_3-y_2x_1,y_4-1+2x_1,y_5-y_4-y_3)$: It is easy to see that at a 
root of $G$, we must have $y_5\!=\!1-2x_1+x^5_1\!=\!0$. The proof 
of Proposition \ref{prop:1} is not much harder: One simply substitutes 
new variables to break down sums with more than $2$ terms and (employing 
the binary expansions of the underlying exponents) monomials 
of degree more than $2$. Proposition \ref{prop:2} follows easily 
upon expanding every complex multiplication (resp.\ complex addition) 
into $4$ real multiplications (resp.\ $2$ real additions), by introducing 
new variables for the real and imaginary parts of the $x_i$. 
 
\medskip 
\noindent 
{\bf Proof of Theorem \ref{thm:cxity}:} First observe that 
$\Log|z|\!\in\!\amoeba(f) \Longleftrightarrow f$ has a complex 
root $\zeta$ with $|\zeta|\!=\!|z|$. Letting $A$ and $B$ denote the 
real and imaginary parts of $f$, and letting $\alpha_i$ and $\beta_i$ denote 
the real and imaginary parts of $\zeta_i$, we thus obtain that 
$\Log|z|\!\in\!\amoeba(f)$ if and only if the polynomial system\\  
\mbox{}\hfill $A(\alpha_1,\beta_1,\ldots, \alpha_n,\beta_n)= 
B(\alpha_1,\beta_1,\ldots,\alpha_n,\beta_n)=0$,  
$\alpha^2_1+\beta^2_1=|z_1|^2,\ldots,\alpha^2_n+\beta^2_n=|z_n|^2$
\hfill \mbox{}\\
has a root 
$(\alpha,\beta)\!=\!(\alpha_1,\ldots,\alpha_n,\beta_1,\ldots,\beta_n)\!\in\!
\R^{2n}$. Now, while the preceding 
system of equations has size significantly larger than 
$\size(z)+\size(f)$ (due to the underlying expansions of powers of 
$\zeta_i\!=\!\alpha_i+\sqrt{-1}\beta_i$), we can introduce new variables 
and equations (via Propositions \ref{prop:1} and \ref{prop:2}) 
to obtain another polynomial system, also with a real solution 
if and only if $\Log|z|\!\in\!\amoeba(f)$, with size linear in 
$\size(z)+\size(f)$ instead. Applying Theorem \ref{thm:canny}, we obtain 
our $\pspa$ upper bound. 

Our $\np$-hardness complexity lower bound follows immediately from  
Theorem \ref{thm:plai}, since $|\zeta|\!=\!1\Longleftrightarrow 
\Log|\zeta|\!=\!0$. \qed   

\begin{rem} 
A reduction of amoeba membership to the Existential Theory 
of the Reals, with an $\expt$ complexity upper bound instead, 
was observed in \cite[Sec.\ 2.2]{theobald}.  
\dia 
\end{rem}

\section{The Proof of Theorem \ref{thm:finally}}
\label{sec:finally} 
We now assume that $f$ is $n$-variate. 
\subsection{Proof of Part (a) of Assertion (1)}  
\label{sub:2} 
When $k\!=\!1$ it is clear that $f$ is of the form $g(x^a)$ for 
some $g\!\in\!\C\!\left[x^{\pm 1}_1\right]$ and 
$a\!\in\!\Z^n\setminus\{\bO\}$. Assertion (1) of Lemma \ref{lemma:basic}, 
and Theorem \ref{thm:harder} applied to $g$, then immediately imply the second 
bound from Part (a). So let us now assume $k\!\geq\!2$. 

Let $w\!:=\!(\log |\zeta_1|,\ldots,\log|\zeta_n|)\!\in\!\amoeba(f)$ and 
assume without loss of generality that 
$|c_1 \zeta^{a_1}|\!\geq\!|c_2\zeta^{a_2}|\!\geq \cdots \geq\!|c_t\zeta^{a_t}|$.
Since $f(\zeta)\!=\!0$ implies that $|c_1\zeta^{a_1}|=|c_2\zeta^{a_2}+\cdots+ 
c_t\zeta^{a_t}|$, the\linebreak 
Triangle Inequality immediately implies 
that $|c_1\zeta^{a_1}|\!\leq\!(t-1)|c_2\zeta^{a_2}|$. 
Taking logarithms, we then obtain 
\begin{eqnarray} 
\label{eqn:3} 
& a_1\cdot w +\log|c_1|\geq \cdots \geq a_t\cdot w +\log|c_t|& \text{ and}  
\end{eqnarray} 
\begin{eqnarray} 
\label{eqn:4} 
& a_1\cdot w +\log|c_1|\leq \log(t-1)+a_2\cdot w +\log|c_2|& 
\end{eqnarray} 
For each $i\!\in\!\{2,\ldots,t\}$ let us then define $\delta_i$ to be
the shortest vector such that \\ 
\mbox{}\hfill 
$a_1\cdot (w+\delta_i)+\log|c_1| 
=a_i\cdot (w+\delta_i)+\log|c_i|$. \hfill\mbox{}\\ 
Note that $\delta_i\!=\!\lambda_i(a_i-a_1)$ for some nonnegative $\lambda_i$ 
since we are trying to affect the dot-product $\delta_i\cdot(a_1-a_i)$. 
In particular, 
$\lambda_i\!=\!\frac{(a_1-a_i)\cdot w+\log|c_1/c_i|}{|a_1-a_i|^2}$ so
that $|\delta_i|\!=\!\frac{(a_1-a_i)\cdot w + \log|c_1/c_i|}
{|a_1-a_i|}$. (Indeed, Inequality (\ref{eqn:3})   
implies that $(a_1-a_i)\cdot w+\log|c_1/c_i|\!\geq\!0$.)  

Inequality (\ref{eqn:4}) implies that 
$(a_1-a_2)\cdot w+\log|c_1/c_2|\!\leq\!\log(t-1)$. 
We thus obtain\linebreak $|\delta_2|\!\leq\!\frac{\log(t-1)}
{|a_1-a_2|}\!\leq\!\log(t-1)$.
So let $i_0\!\in\!\{2,\ldots,t\}$ be any $i$ minimizing
$|\delta_i|$. We of course have $|\delta_{i_0}|\!\leq\!\log(t-1)$, and 
by the definition of $\delta_{i_0}$ we have\\ 
\mbox{}\hfill  
$a_1\cdot (w+\delta_{i_0})+\log|c_1|\!=\!a_{i_0}\cdot (w+\delta_{i_0})
+\log|c_{i_0}|$.\hfill\mbox{}\\ 
Moreover, the fact that $\delta_{i_0}$ is the shortest among the $\delta_i$ 
implies that\\
\mbox{}\hfill $a_1\cdot (w+\delta_{i_0})+\log|c_1|\!\geq\!a_i\cdot (w+ 
\delta_{i_0})+\log|c_i|$\hfill\mbox{}\\  
for all $i$. Otherwise, we would have $a_1\cdot 
(w+\delta_{i_0})+\log|c_1|\!<\!a_i\cdot (w+\delta_{i_0})+\log|c_i|$ and
$a_1\cdot w+\log|c_1|\!\geq\!a_i\cdot w+\log|c_i|$ (the latter 
following from Inequality (\ref{eqn:3})). Taking a convex linear 
combination of the last two inequalities, it is then clear that 
there must be a $\mu\!\in\![0,1)$ such that $a_1\cdot 
(w+\mu\delta_{i_0})+\log|c_1| 
\!=\!a_i\cdot (w+\mu\delta_{i_0})+\log|c_i|$. Thus, by the 
definition of $\delta_i$, we would obtain $|\delta_i|
\!\leq\!\mu|\delta_{i_0}|\!<\!|\delta_{i_0}|$ --- a contradiction.

We thus have the following:\\ 
\mbox{}\hfill $a_1\cdot (w+\delta_{i_0})-(-\log|c_1|)\!=\!
a_{i_0}\cdot (w+\delta_{i_0})-(-\log|c_{i_0}|)$,\hfill \mbox{}\\ 
\mbox{}\hfill $a_1\cdot (w+\delta_{i_0})-(-\log|c_1|)\!\geq\!
a_i\cdot (w+\delta_{i_0})-(-\log|c_i|)$\hfill\mbox{}\\ 
for all $i$, and $|\delta_{i_0}|\!\leq\!\log (t-1)$. 
This implies that $w+\delta_{i_0}\!\in\!\atrop(f)$. In other words, 
we've found a point in $\atrop(f)$ sufficiently near $\Log|\zeta|$ to prove our 
desired upper bound. \qed 

\subsection{Proving Part (b) of Assertion (1)}  
\label{sub:big} 
We begin with a refinement of the special case $n\!=\!1$.
Let $\#S$ denote the cardinality of a set $S$.
\begin{thm}
\label{thm:haus}
Suppose $f$ is any univariate $t$-nomial with $t\!\geq\!3$ and 
$s\!:=\!\#\atrop(f)$. (So $1\!\leq\!s\!\leq\!t-1$.) Then for any
$v\!\in\!\atrop(f)$ there is a root $\zeta\!\in\!\C$ of $f$ with
$|v-\log|\zeta||\!<\!\log 2$, 
$|v-\log|\zeta||\!\leq\!\log\min\{18,t-1\}$, 
or $|v-\log|\zeta||\!<\!(\log 9)s-\log \frac{9}{2}\!<\!2.2s-1.5$, 
according as $s$ is $1$, $2$, 
or $\geq\!2$. In particular, $|v-\log|\zeta||\!<\!(\log 9)t-\log\frac{81}{2}
\!<\!2.2t-3.7$ for all $t\!\geq\!3$. 
\end{thm}

\noindent
{\bf Proof:} Following the notation of Theorem \ref{thm:harder}, let 
$\Gamma$ be the connected component of $U_f$ 
containing $v\!\in\!\atrop(f)$ and $m\!:=\!\#(\Gamma\cap \atrop(f))$.  
(So $1\!\leq\!m\!\leq\!s$.) 
The quantity $|v-\log|\zeta||$ is thus clearly maximized, for instance, when 
$v$ is as far to the left as possible and $\log|\zeta|$ is as far to the 
right as possible. In other words,\\ 
\mbox{}\hfill $|v-\log|\zeta||\!<\!\log(3)+(\log 9)(m-2)
+\log(3)+\delta$,\hfill\mbox{}\\ 
where $\delta$ is $\log 3$ or $\log 2$, according as 
$m\!<\!s$ or $m\!=\!s$. We thus obtain the largest possible upper bound 
of $(\log 9)s-\log\frac{9}{2}$ when $m\!=\!s$. Note also that $s\!\leq\!t-1$. 
So now we merely need to refine the cases with $s\!\in\!\{1,2\}$. 

\scalebox{1}[1]{The case $s\!=\!1$ follows from Corollary 
\ref{cor:extremecauchy} since $\min\atrop(f)\!=\!\max\atrop(f)$ here.} 

The case $s\!=\!2$ proceeds as follows: 
If $m\!=\!1$ then $\Gamma$ is an open interval of width $2\log 3$ 
with $v$ at its median, so we must have $|v-\log|\zeta||\!<\!\log 3$. 
If $m\!=\!2$ then $\Gamma$ is an open interval of width at most $4\log 3$, 
but we still have\\ 
\mbox{}\hfill  
$\min \atrop(f)-\log 2\!<\!\log|\zeta|\!<\!\max \atrop(f)+\log 2$.
\hfill\mbox{}\\   
So $|v-\log|\zeta||$ can again be maximized by having 
$v$ as far left as possible and $\log|\zeta|$ as far right\linebreak 
\scalebox{.935}[1]{as possible. 
In particular, $s\!=\!2$ implies that $\atrop(f)\!=\!\{\min \atrop(f), 
\max\atrop(f)\}$.}\linebreak 
So we obtain $|v-\log|\zeta||\!<\!\log(3)+\log(3)+\log(2)
\!=\!\log 18$. In addition, we can apply Corollary \ref{cor:extremecauchy} 
to observe that there is always a root $\zeta\!\in\!\C$ of $f$ with 
$|\min\atrop(f)-\log|\zeta||\!\leq\!\log(t-1)$, and the same bound 
can be attained for $|\max\atrop(f)-\log|\zeta||$, possibly with a\linebreak  
different root $\zeta$. So we obtain 
$|v-\log|\zeta||\!\leq\!\log\min\{18,t-1\}$. \qed 

\medskip 
We will handle the case $n\!\geq\!2$ by showing that any point 
$v\!\in\!\atrop(f)$ lies close to the intersection of $\amoeba(f)$ with a 
specially chosen line also containing $v$. With some care, this enables us to 
reduce to the case $n\!=\!1$. In particular, intersecting a line with 
$\amoeba(f)$ is the same as evaluating $f$ along a monomial curve, and we'll 
need a technical lemma to pick exponents that permit an easy reduction to 
$n\!=\!1$. 
\begin{thm} 
\label{thm:whew} 
Given any subset $\{a_1,\ldots,a_t\}\!\subset\!\Zn$ of cardinality 
$t\!\geq\!n+1$, there exists an 
$\alpha\!=\!(\alpha_1,\ldots,\alpha_n)\!\in\!\Zn\!\setminus\!\{\bO\}$ 
such that the dot-products 
$\alpha\cdot a_1,\ldots,\alpha\cdot a_t$ are pair-wise distinct 
and, for all $i\!\in\![n]$, $|\alpha_i|\!\leq\!\ceil{\frac{1}{4}t(t-1)}$ 
or $|\alpha_i|\!\leq\!t-2$, according as $n\!\geq\!3$ or $n\!=\!2$.  
\end{thm} 

\noindent 
{\bf Proof:} Observe that for the $\alpha\cdot a_i$ to remain 
distinct we must have $\alpha$  
avoid a set of $\leq\!t(t-1)/2$ hyperplanes, depending on 
$\{a_1,\ldots,a_t\}$. This is equivalent to $\alpha$ avoiding the zero set 
of an $n$-variate polynomial of degree $t(t-1)/2$. Schwartz's Lemma 
(see, e.g., \cite{schwartz}) then tells us that for {\em any} $S\!\subset\!\Z$ 
with $\#S\!>\!t(t-1)/2$ there is an $\alpha\!\in\!S^n$ avoiding
our aforementioned set of hyperplanes. Picking 
$S\!=\!\left\{-\ceil{\frac{1}{4}t(t-1)},\ldots, 
\ceil{\frac{1}{4}t(t-1)}\right\}$ then gives us the 
case $n\!\geq\!3$. 

For the case $n\!=\!2$, it is enough to prove that the set of lattice 
points\\ 
\mbox{}\hfill $X\!:=\!\{-(t-2),\ldots,t-2\}\times \{1,\ldots,t-2\}$\hfill
\mbox{}\\ 
contains at least $1+t(t-1)/2$ distinct directions (and thus we can always 
find a suitable $\alpha\!\in\!X$). In other words, we need to prove  
that $X$ has at least $1+t(t-1)/2$ points with relatively prime coordinates. 
Throwing out the directions $(1,0)$ and $(0,1)$, it is then enough to 
show that $Y\!:=\!\{1,\ldots,t-2\}^2$ contains at least 
$\frac{t(t-1)}{4}-\frac{1}{2}$ points with relatively prime coordinates. 
The number of such points, for arbitrary $t$, forms the sequence 
A018805 in Sloane's Online Encyclopedia of Integer Sequences 
\cite{oeis}. A routine, but tedious calculation then yields the 
$t\!\in\!\{3,\ldots,45\}$ portion of the $n\!=\!2$ case. 

The remaining cases can be settled as follows: By a standard M\"obius 
inversion argument, the number of points with relatively prime 
coordinates in $Y$ is exactly 
$\sum\limits^{t-2}_{d=1}\mu(d)\floor{(t-2)/d}^2$ where $\mu$ is the 
classical {\em M\"obius function} (see, e.g., \cite{hardy}). A simple 
expansion then yields our desired number of points to be bounded from below 
by\\   
\mbox{}\hfill $A(t):=\frac{(t-2)^2}{\zeta(2)}-4(t-2)-2(t-2)\log(t-2)-2\zeta(2)
(t-1)$.\hfill\mbox{}\\ 
A simple derivative calculation then yields that 
$A(t)-\frac{t(t-1)}{4}+\frac{1}{2}$ is increasing for all $t\!\geq\!25$. 
So it's enough to prove that $A(46)\!>\!517$. One can 
check via {\tt Maple} that $A(46)\!>\!519.9$, so we are done. \qed 

\medskip 
\noindent
{\bf Proof of Part (b) of Assertion (1):} 
Let $v\!=\!(v_1,\ldots,v_n)$ be any point of $\atrop(f)$. If 
$v\!\in\!\amoeba(f)$ then there is nothing to prove. So let us assume 
$v\!\not\in\!\amoeba(f)$. Since the case $n\!=\!1$ is
immediate from Proposition \ref{prop:dumb}, and Theorem \ref{thm:haus}  
we will assume henceforth that $n\!\geq\!2$.

So we can reduce to the case $k\!=\!n$, let us temporarily assume that 
$k\!<\!n$. Without loss of generality, we can order 
the variables $x_1,\ldots,x_n$ so that the image of $\newt(f)$ under 
the coordinate projection sending $\Rn$ onto $\R^k\times\{0\}^{n-k}$ 
has dimension $k$ (and the restriction of the projection to 
$\newt(f)$ is a bijection). Define 
$g(x_1,\ldots,x_k)\!:=\!f(x_1,\ldots,x_k,e^{v_{k+1}},\ldots,e^{v_n})$. 
By the definition of $\atrop(f)$, $\max_{i\in[t]} |c_i 
e^{a_i\cdot v}|$ is attained for at least two distinct values of $i$.
By our construction of $g$, this monomial norm condition implies that
$(v_1,\ldots,v_k)\!\in\!\atrop(g)$. Clearly then, if we can find a root 
$(\zeta_1,\ldots,\zeta_k)$ of $g$ with 
$|(v_1,\ldots,v_k)-\Log|(\zeta_1,\ldots,\zeta_k)||\!<\!\eps_{k,t}$, 
then $\zeta\!:=\!(\zeta_1,\ldots,\zeta_k,e^{v_{k+1}},\ldots,e^{v_n})$ will 
be a root of $f$ with $|v-\Log|\zeta||\!<\!\eps_{k,t}$. But finding 
such a $(\zeta_1,\ldots,\zeta_k)$ for $g$ is nothing more than an instance 
of the case where the dimension of the underlying Newton polytope is the 
same as the underlying number of variables. 

\scalebox{.96}[1]{So we may assume $k\!=\!n\!\geq\!2$ henceforth. 
Consider a monomial curve $C(t)\!:=\!(\gamma_1 t^{\alpha_1},
\ldots,\gamma_n t^{\alpha_n})$}\linebreak 
with $\alpha\!=\!(\alpha_1,\ldots,\alpha_n)
\!\neq\!\bO$. (Note that $\{\Log|C(t)|\}_{t\in\C}$ 
is always a line in $\Rn$.) Setting $\gamma_i\!=\!e^{v_i}$ for all $i$ 
we obtain $v\!=\!\Log|C(1)|$, independent of $\alpha$. So let 
us pick $\alpha$ satisfying the conclusion of Theorem \ref{thm:whew} and 
set $h(t)\!:=\!f(C(t))$. Then by the definition of $\atrop(f)$, and 
especially because the $\alpha\cdot a_i$ are pair-wise distinct, we can 
conclude that $h$ has exactly $t$ monomial terms {\em and}  
$0\!\in\!\atrop(h)$. So to find a root $\zeta\!\in\!\Cn$ with $\Log|\zeta|$ 
close to $v$, it's enough to prove that $h$ has a root $\rho$ close to $1$. 
Thanks to Theorem \ref{thm:haus}, we can do the latter, so now we simply 
have to account for metric distortion from specializing $f$ along $C(t)$. 

Taking logarithms, $\amoeba(h)$ containing a point at distance $\eps$ from 
$0$ implies that $\amoeba(f)$ contains a point at distance $\leq\!|\alpha|
\eps$ from $v$. So by the coordinate bounds of Theorem \ref{thm:whew}, we are 
done. \qed 

\subsection{Proof of Assertion (2)}  
\label{sub:ass2} 
To prove Part (a), note that $(1,\ldots,1)/(1-t)$ is a root of $\varphi$ and 
thus $p\!:=\!-\log(t-1)(1,\ldots,1)\!\in\!\amoeba(\varphi)$. Note that 
$\newt(\varphi)$ is 
the standard $n$-simplex $\Delta_n\!\subset\!\Rn$. So, by polyhedral duality 
\cite{ziegler}, and the definition of $\atrop(\varphi)$, we have that 
$\atrop(\varphi)$ is the positive codimension locus of the outer normal fan 
of $\Delta_n$. In particular, 
$\atrop(\varphi)\cap\overline{\R^{t-1}_-}$ is the 
boundary of the negative orthant. So 
the distance from $p$ to $\atrop(\varphi)$ is $\log(t-1)$. 

To prove Part (b), note that $(x_1+1)^{t-k}$ has a unique root of multiplicity 
$t-k$ at\linebreak 
$x_1\!=\!-1$. Recall that the roots of a monic univariate polynomial 
are continuous functions of the coefficients, e.g., \cite[Thm.\ 1.3.1, Pg.\ 
10]{rs}.\footnote{The statement there excludes roots of multiplicity 
equal to the degree of the polynomial, but the proof in fact works 
in this case as well.} So then, for any $\eps\!>\!0$,  
we can find a $\delta_\eps\!>\!0$ so that for all $\delta\!\in\!\C$ 
with $|\delta|\!\in\![0,\delta_\eps)$, {\em all} the roots $\zeta_1$ of 
$(x_1+1)^{t-k}-\delta$ satisfy $|\zeta_1+1|\!<\!\eps$. Clearly then, for 
any $\eps'\!>\!0$, taking 
$|\rho_2|,\ldots,|\rho_n|$ sufficiently 
small (or $u_2\!:=\!\log|\rho_2|,\ldots,u_n\!:=\!\log|\rho_n|$ sufficiently 
negative) implies that the distance from any point 
$u\!\in\!\amoeba(f)$ of the form $(u_1,u_2,\ldots,u_n)$ to the hyperplane 
$\{0\}\times \R^{n-1}$ is at most $\eps'$: Simply take $\eps$ so that 
$\eps'\!=\!\log(1+\eps)$ and $|x_2|+ \cdots+|x_n|\!<\!\delta_\eps$.  

On the other hand, by the log-concavity of the binomial coefficients, 
$\anewt\!\left((x_1+1)^{t-k}\right)$ must have an edge of slope $t-k$. 
This will enable us to prove that $\atrop(\psi)$ contains a ray of the 
form $\{(\log(t-k),N,\ldots,N)\}_{N\rightarrow+\infty}$. 
and thus conclude: The points along this ray have distance to $\amoeba(\psi)$ 
approaching $\log(t-k)$, by the preceding paragraph. 

To see why such a ray lies in $\atrop\!\left((x_1+1)^{t-k}\right)$  
simply note that as $N\longrightarrow -\infty$, the linear form 
$\log(t-k)u_1+Nu_2+\cdots+Nu_n-u_{n+1}$ is maximized exactly at the vertices\\
\mbox{}\hfill  $\left(t-k-1,0,\ldots,0,-\log(t-k)\right)$ 
\ \ and \ \ $\left(t-k,0,\ldots,0,0\right)$\hfill\mbox{}\\  
of $\anewt\!\left((x_1+1)^{t-k}\right)$. 
(Indeed, the only other possible vertices of 
$\anewt\!\left((x_1+1)^{t-k}\right)$ are the basis vectors 
$e_2,\ldots,e_k$ of $\R^{n+1}$.) So, by Lemma \ref{lemma:trop}, we are 
done. \qed

\section{Proving Theorem \ref{thm:cxity2}} 
\label{sec:cxity2} 
Let us first recall the following result on comparing monomials in rational 
numbers. 
\begin{thm} \cite[Sec.\ 2.4]{thresh} 
\label{thm:linlog} 
Suppose $\alpha_1,\ldots,\alpha_N\!\in\!\Q$ are positive and 
$\beta_1,\ldots,\beta_N\!\in\!\Z$. Also let $A$ be the  
maximum of the numerators and denominators of the $\alpha_i$ (when 
written in lowest terms) and $B\!:=\!\max_i\{|\beta_i|\}$. Then, within\\
\mbox{}\hfill $O\!\left(N30^N\log(B)(\log \log B)^2\log\log\log(B)
(\log(A)(\log \log A)^2\log\log\log A)^N\right)$\hfill \mbox{} \\
bit operations, we can determine the sign of $\alpha^{\beta_1}_1\cdots 
\alpha^{\beta_N}_N-1$. \qed
\end{thm} 

\noindent 
While the underlying algorithm is a simple application of 
Arithmetic-Geometric\linebreak Mean Iteration (see, e.g., 
\cite{dan}), its complexity bound hinges on a deep estimate 
of Nesterenko \cite{nesterenko}, which in turn refines seminal  
work of Matveev \cite{matveev} and Alan Baker \cite{baker} on 
linear forms in logarithms.  

\medskip 
\noindent 
{\bf Proof of Theorem \ref{thm:cxity2}:} 
From Lemma \ref{lemma:trop}, it is clear that we 
merely need an efficient method to compare quantities of the 
form $|c_iz^{a_i}|$, and there are exactly $t-1$ such comparisons 
to be done. So our first complexity bound follows immediately 
from the special case of Theorem \ref{thm:linlog} where 
$A\!=\!2^\sigma$, $B\!=\!d$, and $N\!=\!2n+2$. 
In particular, $30\log 2\!<\!20.8$. 

The second assertion follows almost trivially: Thanks to the 
exponential form of the\linebreak 
coefficients and the query point, one can take logarithms to 
reduce to comparing integer linear combinations of rational numbers of bit 
size linear in $\max\{\sigma,\log d\}$. So the underlying monomial norm 
comparisons can 
be reduced to standard techniques 
for fast integer multiplication (see, e.g., \cite[Pg.\ 43]{bs}). \qed 

\section*{Acknowledgements} 
We thank Pascal Koiran and Natacha Portier at the \'Ecole Normale Sup\'erieure 
of Lyon, and Peter Gritzmann at Technical University of M\"unchen, 
for their splendid hospitality while this paper was completed. Rojas also 
gratefully acknowledge the support of Labex MILYON and the Bavarian Ministry 
of Science, Research, and Art.

We also thank Eleanor Anthony, Frederic Bihan, Sheridan Grant, and 
Timo de Wolff for comments on earlier drafts of this paper. Special thanks go 
to Jan-Erik Bj\"ork and Jean-Yves Welschinger for respectively bringing 
Ostrowski's paper \cite{ostrowski} and Viro's paper 
\cite{virologpaper} to our attention, to Matt Young for pointing out 
the relevance of M\"obius inversion for counting lattice directions, 
and to Jens Forsg\aa{}rd for pointing out Pellet's Theorem. 

As this paper was being developed, we heard the tragic news that 
Mikael Passare died in a hiking accident in Oman on September 15, 2011. 
Mikael was not only a titan in the study of amoebae. He brightened 
the lives of everyone who knew him. He is sorely missed. 

\small 
\bibliographystyle{acm}

\end{document}